\documentclass{article}
\usepackage{amsfonts}

\usepackage{graphicx}
\usepackage{amsmath}
\usepackage[symbol]{footmisc}
\usepackage[french]{babel}

\newtheorem{theorem}{Theorem}[section]

\newtheorem{corollary}[theorem]{Corollary}

\newtheorem{definition}[theorem]{Definition}

\newtheorem{lemma}[theorem]{Lemma}

\newtheorem{proposition}[theorem]{Proposition}
\newtheorem{remark}[theorem]{Remark}

\newenvironment{proof}[1][Proof]{\textbf{#1.} }{\ \rule{0.5em}{0.5em}}
\input{tcilatex}

\begin{document}

\begin{center}
{\huge q-power symmetric functions and q-exponential formula}

\bigskip

{\large Vincent Brugidou}$^{a}$ 
\footnotetext{$^{a}$\textit{\ }vincent.brugidou@univ-lille.fr}
\end{center}

\bigskip

\begin{center}
\textit{Univ. Lille, CNRS, UMR 8524 - Laboratoire Paul Painlev\'{e}, F-
59000 Lille, France}
\end{center}

\bigskip

\textbf{Abstract: }Let $\lambda =\left( \lambda _{1},\lambda
_{2},...,\lambda _{r}\right) $ be an integer partition, and $\left[
p_{\lambda }\right] $ the $q$-analog of the symmetric power function $%
p_{\lambda }$. This $q$-analogue has been defined as a special case, in the
author's previous article: ''A $q$-analog of certain symmetric functions and
one of its specializations''. Here, we prove that a large part of the
classical relations between $p_{\lambda }$, on one hand, and the elementary\
and complete symmetric functions $e_{n}$ and $h_{n}$, on the other hand,
have $q$-analogues with $\left[ p_{\lambda }\right] $. In particular, the
generating functions $E\left( t\right) =\sum\nolimits_{n\geq 0}e_{n}t^{n}$
and $H\left( t\right) =\sum\nolimits_{n\geq 0}h_{n}t^{n}$ are expressed in
terms of $\left[ p_{n}\right] $,\ using Gessel's $q$-exponential formula and
a variant of it. A factorization of these generating functions into infinite 
$q$-products, which has no classical counterpart, is established. By
specializing these results, we show that the $q$-binomial theorem is a
special case of these infinite $q$-products. We also obtain new formulas for
the tree inversions enumerators and for certain $q$-orthogonal polynomials,
detailing the case of dicrete $q$-Hermite polynomials.

\bigskip

\textit{keywords: }Symmetric functions, $q$-calculus, exponential formula,
tree inversions, $q$-orthogonal polynomials.

\section{\protect\bigskip \textbf{Introduction}}

Let $p_{n}^{\left( r\right) }$\ be the symmetric functions, defined for any
pair of integers $\left( n,r\right) $ such that $n\geq r\geq 0$, \ by $%
p_{n}^{\left( r\right) }=\delta _{n}^{r}$ (Kronecker symbol) if $r=0$, and
if $r>0$, by

\begin{equation}
p_{n}^{(r)}=\sum\limits_{\left| \lambda \right| =n,\;l\left( \lambda \right)
=r}m_{\lambda }\text{,}  \tag{1.1}
\end{equation}

where $m_{\lambda }$ represents symmetric monomial functions, the sum being
over the integer partitions $\lambda $ of $n$, with length $l(\lambda )=r$.
These functions $p_{n}^{(r)}$ are introduced with this notation in [13,
Exercise 19, p. 33]. In [4], we defined their q-analogs $\left[
p_{n}^{\left( r\right) }\right] $ for any pair of integers $\left(
n,r\right) $ such that $n\geq r\geq 0$, and proved some properties of these $%
q$-analogs. In the special case $r=1$, $\left[ p_{n}^{\left( 1\right) }%
\right] $ is the $q$-analog of the power symmetric function $p_{n}$ and also
denoted as $\left[ p_{n}\right] $. Additionally, we have also defined, for
any integer partition $\lambda =\left( \lambda _{1},\lambda _{2},...,\lambda
_{r}\right) $, the $q$-analog of $p_{\lambda }$ as follows:

\begin{equation}
\left[ p_{\lambda }\right] =\left[ p_{\lambda _{1}}\right] \left[ p_{\lambda
_{2}}\right] ...\left[ p_{\lambda _{r}}\right]  \tag{1.2}
\end{equation}

The first objective of this article is to complete [4] by showing that a
large part of the classical relations linking power symmetric functions to
elementary and complete symmetric functions, or to their generating series $%
E(t)$ and $H(t)$, have $q$-analogs with $\left[ p_{n}\right] $ or $\left[
p_{\lambda }\right] .$ The second objective is to show that $E(t)$ and $H(t)$
can be expressed in terms of $\left[ p_{n}\right] $\ as infinite $q$%
-products without classical counterparts. A final objective is to give some
applications of these results by specialization.

The organization of the article is as follows. In section 2, we establish
the notations and recall the necessary prerequisites for symmetric
functions, integer partitions, and $q$-calculus. In Section 3, we recall the
definitions of $\left[ p_{n}^{\left( r\right) }\right] $ and the relations
between $\left[ p_{n}\right] $ and $e_{n}$ given in [4], then we establish
the relations between $\left[ p_{n}\right] $ and $h_{n}$. In Section 4, we
obtain q-analogs of the following classical formulas (see, for example, $%
\left[ \text{13, Chap.1, Eq. (2.14)}\right] $) 
\begin{equation}
e_{n}=\sum\limits_{\left| \lambda \right| =n}\varepsilon _{n}z_{\lambda
}^{-1}p_{\lambda }\text{,}  \tag{1.3}
\end{equation}

\begin{equation}
h_{n}=\sum\limits_{\left| \lambda \right| =n}z_{\lambda }^{-1}p_{\lambda }%
\text{,}  \tag{1.4}
\end{equation}
by introducing $q$-deformations of $z_{\lambda }$, the cardinality of the
centralizer for permutations of type $\lambda .$

In sections 5 and 6, we obtain q-analogs of classical expressions in which
the generating series $E(t)$ and $H(t)$ are formal compositions of the
exponential series with the series $-P(-t)$ and $P(t))$, respectively, where

\begin{equation}
P(t)=\sum\limits_{n\geq 1}\dfrac{p_{n}}{n}t^{n}\text{.}  \tag{1.5}
\end{equation}

In the case of $E(t)$, we use the $q$-composition introduced by Gessel in
[8] for that purpose. For $H(t)$, we introduce another $q$-composition, for
which we show a combinatorial interpretation similar to that presented in $%
\left[ 8\right] $. We thus obtain a variant of Gessel's $q$-exponential
formula, in the sense that like this one, it is a $q$-analog of the
classical exponential formula (for the latter, see [15, chap. 5]). In
Section 7, we highlight the functionnal link between these two $q$%
-exponential formulas.

In Section 8, we show that $E(t)$ and $H(t)$ can be expressed in terms of $%
\left[ p_{n}\right] $\ as infinite $q$-products.These products have no
classical counterpart, i.e., when $q=1$. Section 9 gives examples of how to
apply the formulas obtained in the previous sections to various $q$-series,
by considering them as specializations in the sense of symmetric functions.
Thus, we show that the $q$-binomial theorem is a special case of the above
infinite $q$-product for $E(t)$. By identifying the generating functions of
\ certain $q$-orthogonal polynomials given in [11] with the infinite $q$%
-products mentionned above, we also derive new formulas for these
polynomials. An example of formulas obtained for the discrete $q$-Hermite I
polynomials, is for $n\geq 2$ :

\begin{equation}
\text{\ }H_{n}\left( x;q\right) =\left| 
\begin{array}{ccccccc}
x & -\left( 1-q\right) & 0 & 0 & . & . & 0 \\ 
x^{2}-1 & qx & -(1-q^{2}) & 0 & . & . & 0 \\ 
x\left( x^{2}-1\right) & q\left( x^{2}-1\right) & q^{2}x & -\left(
1-q^{3}\right) & . & . & . \\ 
. & . & . & . & . & . & . \\ 
. & . & . &  & . & . & 0 \\ 
. & . & . &  &  &  & -\left( 1-q^{n-1}\right) \\ 
x^{n-2}\left( x^{2}-1\right) & qx^{n-3}\left( x^{2}-1\right) & 
q^{2}x^{n-4}\left( x^{2}-1\right) & . & . & . & q^{n-1}x
\end{array}
\right|  \tag{1.6}
\end{equation}

Various potential developments are indicated throughout the article and in
Section 10.

\section{Prerequisites}

\QTP{Body Math}
First some standard notations are recalled. $\delta _{i}^{j}$ is the
Kronecker symbol equal to 1 if $i=j$ and $0$ otherwise. If $E$ is a finite
set, $\left| E\right| $ denotes its cardinality. $\mathbb{N},\mathbb{Z,Q,C}$
denote natural integers, integers, rational numbers and complex numbers,
respectively. If $A$ is a commutative ring (with identity), $A[[t]]$ and $A%
\left[ \left[ t_{1},t_{2,...}\right] \right] $ are the rings of formal power
series (abbreviated as fps) with coefficients in $A$, respectively in the
indeterminate $t$ and in the (finite or infinite) sequence of indeterminates 
$t_{1},t_{2},...$. The order (or valuation) of a fps $F\neq 0$ is the
smallest integer $n$, denoted $\omega \left( F\right) $, such that the
homogeneous component of degree $n$ of $F$ is nonzero. The convergences
considered in this article are always those of the fps for the valuation
topology (see for example [3, Chap.4 \S 4] or [6, \S 1.12]). We will refer
to it as formal convergence. For $n\in \mathbb{N}^{\ast }=\mathbb{%
N\backslash }\left\{ 0\right\} $, we set $\mathbf{n=}\left\{
1,2,...,n\right\} $.

\QTP{Body Math}
For the integer partitions and symmetric functions, the notations are those
of [13, Chap.1] or [15, Chap.7], to which we refer the reader for further
details. If $\lambda =\left( \lambda _{1},\lambda _{2},...,\lambda
_{r}\right) $ is a partition of $n\in \mathbb{N}$, let $\left| \lambda
\right| =\lambda _{1}+\lambda _{2}+...+\lambda _{r}=n,\;l(\lambda )=r$ and$%
\; $

\QTP{Body Math}
\begin{equation}
n(\lambda )=\sum\limits_{i=1}^{r}(i-1)\lambda _{i}=\sum\limits_{i\geq 1}%
\binom{\lambda _{i}^{\prime }}{2}  \tag{2.1}
\end{equation}
where $\lambda ^{\prime }=\left( \lambda _{1}^{\prime },\lambda _{2}^{\prime
},..\right) $ is the conjugate of the partition $\lambda $. Note that we
have $\left( \lambda ^{\prime }\right) ^{\prime }=\lambda $ and then $%
n(\lambda ^{\prime })=\sum\nolimits_{i\geq 1}(i-1)\lambda _{i}^{\prime
}=\sum\nolimits_{i=1}^{r}\binom{\lambda _{i}}{2}$

\QTP{Body Math}
We generalize Eq. $\left( 2.1\right) $ to any strict composition of $n$,
i.e. any r-multiplet $u=\left( u_{1},u_{2},...,u_{r}\right) $ of integers
strictly greater than zero such that $u_{1}+u_{2}+...+u_{r}=n$. This is
achieved by defining $\left| u\right| =u_{1}+u_{2}+...+u_{r}$, $l(u)=r$ and 
\begin{equation}
n\left( u\right) =\sum\limits_{i\geq 1}(i-1)u_{i}  \tag{2.2}
\end{equation}

\QTP{Body Math}
$L$et's now return to integer partitions. For $\lambda =\left( \lambda
_{1},\lambda _{2},...,\lambda _{r}\right) $ and $i\in \mathbb{N}^{\ast }$, $%
m_{i}\left( \lambda \right) $ represents the number of parts of $\lambda $
equal to $i$, $m\left( \lambda \right) =\left( m_{i}\left( \lambda \right)
\right) _{i\geq 1}$ is the sequence of multiplicities of $\lambda $ (almost
all zero), and $z_{\lambda }$ is given by:

\QTP{Body Math}
\begin{equation}
z_{\lambda }=\lambda _{1}\lambda _{2}...\lambda _{r}\prod\limits_{i\geq
1}\left( m_{i}\left( \lambda \right) \right) !=\prod\limits_{i\geq
1}i^{m_{i}\left( \lambda \right) }(m_{i}\left( \lambda \right) !)  \tag{2.3}
\end{equation}

\QTP{Body Math}
If $K$ is a commutative field, $\Lambda _{K}$ is the algebra of symmetric
functions in indeterminates $X=\left( x_{i}\right) _{i\geq 1}$ with
coefficients in $K$. For the purpose of this article, $K$ will be a field
containing $\mathbb{Q}$, typically $\mathbb{Q(}q)$ . For $\lambda $
describing the set of partitions, $\left( m_{\lambda }\right) ,\left(
e_{\lambda }\right) ,\left( h_{\lambda }\right) ,\left( p_{\lambda }\right) $%
, are the four classical bases of $\Lambda _{K}$. By agreeing that $%
e_{0}=h_{0}=1\;$we recall that in $\Lambda _{K}\left[ \left[ t\right] \right]
$

\QTP{Body Math}
\begin{equation}
E(t)=\sum\limits_{n\geq 0}e_{n}t^{n}=\prod\limits_{i\geq
1}(1+x_{i}t)\;;\;\;H(t)=\sum\limits_{n\geq 0}h_{n}t^{n}=\left( E(-t)\right)
^{-1}=\dfrac{1}{\prod\limits_{i\geq 1}\left( 1-x_{i}t\right) }  \tag{2.4}
\end{equation}

\QTP{Body Math}
For the $q$-analogs we use the notations of [10]. The index $q$ can be
omitted if there \ is no ambiguity. For ($n,k)\in \mathbb{N}^{2}$, we have $%
\left[ n\right] _{q}=1+q+q^{2}+...+q^{n-1}=1-q^{n}/1-q$ for $n\neq 0\;$and $%
\left[ 0\right] _{q}=0$; $\left[ n\right] _{q}!=\left[ 1\right] \left[ 2%
\right] ...\left[ n\right] $ for $n\neq 0\;$and $\left[ 0\right] _{q}!=1.$%
Then

\QTP{Body Math}
\begin{equation*}
\QATOPD[ ] {n}{k}_{q}=\dfrac{\left[ n\right] !}{\left[ k\right] !\left[ n-k%
\right] !}\text{ for }n\geq k\geq 0\text{ \ \ and }\QATOPD[ ] {n}{k}_{q}=0\;%
\text{otherwise.}
\end{equation*}
$D^{r}F(t)$ is the usual derivation of order $r$ with respect to $t$ of the
fps $F(t).$The $q$-d\'{e}rivation of $F(t)$ is defined by: 
\begin{equation*}
D_{q}F(t)=\dfrac{F(qt)-F(t)}{(q-1)t}\text{ \ \ and for }r\geq
1\;\;D_{q}^{r}=D_{q}\left( D^{r-1}\right) \text{ \ with \ }D_{q}^{0}F(t)=F
\end{equation*}

\QTP{Body Math}
In particular, $D_{q}^{r}t^{n}=\left[ r\right] !\QATOPD[ ] {n}{r}t^{n-r}$
for $n\geq r$ and $D_{q}t^{0}=0$.

\QTP{Body Math}
We also have the following formula for the quotient (see [10, p.3]), where $%
A(t)$ and $B(t)$ are fps with $B(0)\neq 0$ 
\begin{equation}
D_{q}\left( \dfrac{A(t)}{B(t)}\right) =\dfrac{B(t)D_{q}A(t)-A(qt)D_{q}B(t)}{%
B(t)B(qt)}  \tag{2.5}
\end{equation}

\QTP{Body Math}
$\bigskip $

\section{$\protect\bigskip $Relationships between $\left[ p_{n}\right] $
and, $e_{n}$ or $h_{n}$}

\QTP{Body Math}
$\bigskip $

\subsubsection{$\protect\bigskip $Reminder on $\left[ p_{n}^{\left( r\right)
}\right] $}

\QTP{Body Math}
For the reader's convenience, we recall the definition of $\left[ p_{n}^{(r)}%
\right] $ as introduced in [4]. For any pair of integers $\left( n,r\right) $%
\ such that $n\geq r\geq 0$, we defined $\left[ p_{n}^{(r)}\right] $ by: 
\begin{equation*}
\sum\limits_{n\geq r}\left[ p_{n}^{(r)}\right] _{q}\left( -t\right) ^{n-r}=%
\dfrac{1}{\left[ r\right] _{q}!}\dfrac{D_{q}^{r}E(t)}{E(t)}.
\end{equation*}
In particular, this definition implies $\left[ p_{n}^{(0)}\right]
_{q}=\delta _{n}^{0}$.

We have deduced the system of equations (the subscript q is implied): 
\begin{equation*}
\text{For }n\geq r\;\sum\nolimits_{k=r}^{n}\left( -1\right) ^{k-r}e_{n-k}%
\left[ p_{k}^{(r)}\right] =\QATOPD[ ] {n}{r}e_{n}\text{,}
\end{equation*}
whose solution is given by the following proposition:

\begin{proposition}
For $n\geq r\geq 1$ we have: 
\begin{equation*}
\left[ p_{n}^{(r)}\right] =\left| 
\begin{array}{cccccccc}
\QATOPD[ ] {r}{r}e_{r} & 1 & 0 & 0 & . & . & . & 0 \\ 
\QATOPD[ ] {r+1}{r}e_{r+1} & e_{1} & 1 & 0 & . & . & . & 0 \\ 
\QATOPD[ ] {r+2}{r}e_{r+2} & e_{2} & e_{1} & 1 & 0 & . & . & 0 \\ 
. & . &  & . & . & . & . & . \\ 
. & . &  &  & . & . & . & . \\ 
. & . &  &  &  & . & . & 0 \\ 
. & . &  &  &  &  & . & 1 \\ 
\QATOPD[ ] {n}{r}e_{n} & e_{n-r} & e_{n-r-1} & . & . & . & e_{2} & e_{1_{{}}}
\end{array}
\right| 
\end{equation*}
\end{proposition}

This determinant can serve as an alternative definition of $\left[
p_{n}^{(r)}\right] $.

\subsection{Relationships between $\left[ p_{n}\right] $ and $e_{n}$
(reminder)}

\QTP{Body Math}
When\textbf{\ } $r=1$, the previous formulas become, with $\left[ p_{n}^{(1)}%
\right] =\left[ p_{n}\right] $, where $p_{n}$ is the power symmetric
function:

\begin{equation}
\sum\limits_{n\geq 1}\left[ p_{n}\right] _{q}\left( -t\right) ^{n-1}=\dfrac{%
D_{q}E(t)}{E(t)},  \tag{3.1}
\end{equation}

and for the system of equations:

\begin{equation}
\text{For }n\geq 1\text{ }\;\;\sum\limits_{k=1}^{n}\left( -1\right)
^{k-1}e_{n-k}\left[ p_{k}\right] =\left[ n\right] e_{n}  \tag{3.2}
\end{equation}

which conversely yields the $q$-Girard-Newton's identities:

\begin{equation}
\text{For }n\geq 1\text{ }\;\;\left[ p_{n}\right] =\sum\limits_{k=1}^{n-1}%
\left( -1\right) ^{k-1}e_{k}\left[ p_{n-k}\right] -\left( -1\right) ^{n-1}%
\left[ n\right] e_{n}\text{.}  \tag{3.3}
\end{equation}

From Eq. $\left( 3.2\right) $ and $\left( 3.3\right) $, we obtained the
following determinant formulas, using Cramer's rule:

\begin{equation}
\left[ p_{n}\right] =\left| 
\begin{array}{ccccccc}
\left[ 1\right] e_{1} & 1 & 0 & 0 & . & . & 0 \\ 
\left[ 2\right] e_{2} & e_{1} & 1 & 0 & . & . & 0 \\ 
. & . & . & . & . & . & . \\ 
. & . & . & . & . & . & . \\ 
. & . & . &  & . & . & 0 \\ 
. & . & . &  &  & . & 1 \\ 
\left[ n\right] e_{n} & e_{n-1} & e_{n-2} & . & . & e_{2} & e_{1}
\end{array}
\right|  \tag{3.4}
\end{equation}

\begin{equation}
\left[ n\right] !e_{n}=\left| 
\begin{array}{ccccccc}
\left[ p_{1}\right] & \left[ 1\right] & 0 & 0 & . & . & 0 \\ 
\left[ p_{2}\right] & \left[ p_{1}\right] & \left[ 2\right] & 0 & . & . & 0
\\ 
. & . & . & . & . & . & . \\ 
. & . & . & . & . & . & . \\ 
. & . &  & . & . & . & 0 \\ 
. & . &  &  &  & . & \left[ n-1\right] \\ 
\left[ p_{n}\right] & \left[ p_{n-1}\right] & \left[ p_{n-2}\right] & . & .
& \left[ p_{2}\right] & \left[ p_{1}\right]
\end{array}
\right|  \tag{3.5}
\end{equation}

We observe that Eq. $\left( 3.4\right) $ and $\left( 3.5\right) $\ are
q-analogs of the classical case; see for instance the first two equations of
[13, Chap. 1, Example 8\ p.28].

\subsection{Relationships between $\left[ p_{n}\right] $ and $h_{n}$}

Let $A(t)=1,$ $B(t)=F(t)$, $G(t)=\left( F(t)\right) ^{-1}$, and suppose that 
$F(0)\neq 0$ which implies that $F$ and $G$ are invertible. Eq. $\left(
2.5\right) $ gives 
\begin{equation*}
\dfrac{D_{q}F(t)}{F(t)}=-\dfrac{D_{q}G(t)}{G(qt)}\text{.}
\end{equation*}
If we choose for $F(t),$ the generating series $E(t)$ defined by $\left(
2.4\right) $, then $G(t)=H\left( -t\right) $. Therefore $%
D_{q}E(t)/E(t)=-D_{q}\left( H\left( -t\right) \right) /H\left( -qt\right) $.
This, combined with $\left( 3.1\right) $, results in: \ \ 
\begin{equation}
\sum\limits_{n\geq 1}\left[ p_{n}\right] t^{n-1}=\dfrac{D_{q}H(t)}{H(qt)}%
\text{.}  \tag{3.6}
\end{equation}

\bigskip Consequently, 
\begin{equation*}
\left( \sum\nolimits_{n\geq 0}h_{n}\left( qt\right) ^{n}\right) \left(
\sum\nolimits_{n\geq r}\left[ p_{n}^{\left( r\right) }\right] t^{n-1}\right)
=\sum\nolimits_{n\geq r}\left[ n\right] h_{n}t^{n-1}\text{.}
\end{equation*}
By performing the above Cauchy product, we obtain by equating the
coefficients of $t^{n}$, the system of equations :

\begin{equation}
\text{For }n\geq 1\;\;\;\sum\limits_{k=1}^{n}h_{n-k}\left[ p_{k}\right]
q^{n-k}=\left[ n\right] h_{n}\text{.}  \tag{3.7}
\end{equation}

\bigskip From these equations, determinant formulas follow by Cramer's rule:

\begin{equation}
\left( -1\right) ^{n-1}\left[ p_{n}\right] =\left| 
\begin{array}{ccccccc}
\left[ 1\right] h_{1} & 1 & 0 & 0 & . & . & 0 \\ 
\left[ 2\right] h_{2} & h_{1}q & 1 & 0 & . & . & 0 \\ 
\left[ 3\right] h_{3} & h_{2}q^{2} & h_{1}q & 1 & 0 & . & . \\ 
. & . & . & . & . & . & . \\ 
. & . & . &  & . & . & 0 \\ 
. & . & . &  &  & . & 1 \\ 
\left[ n\right] h_{n} & h_{n-1}q^{n-1} & h_{n-2}q^{n-2} & . & . & h_{2}q^{2}
& h_{1}q
\end{array}
\right|  \tag{3.8}
\end{equation}
By inverting $\left( 3.7\right) $, we obtain the system of equations \ 
\begin{equation*}
\text{For }n\geq 1\;\;\left[ p_{n}\right] =\sum\limits_{k=1}^{n-1}-h_{k}%
\left[ p_{n-k}\right] q^{k}+\left[ n\right] h_{n}\text{,}
\end{equation*}
which, according to Cramer's rule, gives:

\begin{equation}
\left[ n\right] !h_{n}=\left| 
\begin{array}{ccccccc}
\left[ p_{1}\right] & -\left[ 1\right] & 0 & 0 & . & . & 0 \\ 
\left[ p_{2}\right] & q\left[ p_{1}\right] & -\left[ 2\right] & 0 & . & . & 0
\\ 
\left[ p_{3}\right] & q\left[ p_{2}\right] & q^{2}\left[ p_{1}\right] & - 
\left[ 3\right] & 0 & . & . \\ 
. & . & . & . & . & . & . \\ 
. & . &  & . & . & . & 0 \\ 
. & . &  &  &  & . & -\left[ n-1\right] \\ 
\left[ p_{n}\right] & q\left[ p_{n-1}\right] & q^{2}\left[ p_{n-2}\right] & .
& . & q^{n-2}\left[ p_{2}\right] & q^{n-1}\left[ p_{1}\right]
\end{array}
\right|  \tag{3.9}
\end{equation}

\bigskip We observe that Eq. $\left( 3.8\right) $ and $\left( 3.9\right) $\
are q-analogs of the classical case; see for instance the last two equations
of [13, Chap.1, Example 8 p.28].

\section{\protect\bigskip Expression of $e_{n}$ and $h_{n}$ as a function of 
$\left[ p_{\protect\lambda }\right] $\ }

We are looking for the q-analog of the classical formulas $\left( 1.3\right) 
$ and $\left( 1.4\right) $, where $\varepsilon _{\lambda }=\left( -1\right)
^{\left| \lambda \right| -l(\lambda )}$ and $z_{\lambda }$ is given by $%
\left( 2.3\right) $. Equivalently, we aim to find $q-$analogs of the
following classical formulas, see for instance [13, \ Chap.1, Eq. $\left(
2.14\right) $ p.25]: 
\begin{equation}
E(t)=\sum\limits_{\lambda }\varepsilon _{\lambda }\;z_{\lambda
}^{-1}p_{\lambda }t^{\left| \lambda \right| }\text{.}  \tag{4.1}
\end{equation}

\begin{equation}
H(t)=\sum\limits_{\lambda }z_{\lambda }^{-1}p_{\lambda }t^{\left| \lambda
\right| }\text{.}  \tag{4.2}
\end{equation}

\subsection{Case of $e_{n}$}

\bigskip We start by proving the following lemma.

\begin{lemma}
For $n\geq 1$ we have 
\begin{equation}
\left[ n\right] !e_{n}=\sum\limits_{k}\left( -1\right) ^{n-r}\dfrac{\left[
n-1\right] !}{\left[ k_{1}\right] \left[ k_{2}\right] ...\left[ k_{r-1}%
\right] \left( \left[ k_{r}\right] !\right) }\left[ p_{n-k_{1}}\right] \left[
p_{k_{1}-k_{2}}\right] ...\left[ p_{k_{r-1}-k_{r}}\right]   \tag{4.3}
\end{equation}
where the sum extends over all sequences of integers $k=\left(
k_{1},k_{2},...,k_{r}\right) $ with $r\in \mathbb{N}^{\ast }$ and $n-1\geq
k_{1}>k_{2}>...>k_{r-1}>k_{r}=0$
\end{lemma}

\begin{remark}
In Eq. $\left( 4.3\right) $, although $k_{r}$ seems unnecessary, we need it
to account for the case where the sequence $\left(
k_{1},k_{2},...,k_{r}\right) $ reduces to $k_{1}=0$.
\end{remark}

\begin{proof}
By induction. This is true for $n=1$. Indeed, $\left[ 1\right] !e_{1}=p_{1}$%
, and the set of $k$ in the sum of $\left( 4.3\right) $ is reduced to $r=1$
and $k_{1}=0$, so this sum is also $(0!/0!)\left[ p_{1}\right] =p_{1}$.
Assuming true $\left( 4.3\right) $ for $n\geq j\geq 1$, let's prove it for $%
n+1$. The $(n+1)$th equation $\left( 3.2\right) $ is $\left[ n+1\right]
e_{n+1}=e_{n}\left[ p_{1}\right] -e_{n-1}\left[ p_{2}\right] +e_{n-2}\left[
p_{3}\right] -...+\left( -1\right) ^{n}\left[ p_{n+1}\right] $. Multiplying
by $\left[ n\right] !$, this gives 
\begin{equation}
\left[ n+1\right] !e_{n+1}=\sum\limits_{j=0}^{n}\left( -1\right) ^{n-j}%
\dfrac{\left[ n\right] !}{\left[ j\right] !}\left[ p_{n+1-j}\right] \left[ j%
\right] !e_{j}\text{.}  \tag{4.4}
\end{equation}
Hence, by the induction hypothesis on $e_{j}$, $j\leq n$:

\begin{equation*}
\left[ n+1\right] !e_{n+1}=\sum\limits_{j=0}^{n}\left( -1\right) ^{n-j}%
\dfrac{\left[ n\right] !}{\left[ j\right] !}\left[ p_{n+1-j}\right]
\sum\limits_{l}\left( -1\right) ^{j-s}\dfrac{\left[ j-1\right] !}{\left[
l_{1}\right] \left[ l_{2}\right] ...\left[ l_{s-1}\right] \left( \left[ l_{s}%
\right] !\right) }\left[ p_{j-l_{1}}\right] \left[ p_{l_{1}-l_{2}}\right] ...%
\left[ p_{l_{s-1}-l_{s}}\right] \text{,}
\end{equation*}
where the second sum is over the set of sequences of integers $l=\left(
l_{1},l_{2},...,l_{s}\right) $, with $j-1\geq l_{1}>l_{2}>...>l_{s}=0$. By
replacing the index $l$\ by $k=\left( k_{1},k_{2},...,k_{r}\right) $, with $%
r=s+1$,\ and $k_{1}=j$, $k_{2}=l_{1}$,$...$, $k_{r}=l_{s}$, we obtain: 
\begin{equation*}
\left[ n+1\right] !e_{n+1}=\sum\limits_{k}\left( -1\right) ^{n+1-r}\dfrac{%
\left[ n\right] !}{\left[ k_{1}\right] \left[ k_{2}\right] ...\left[ k_{r-1}%
\right] \left( \left[ k_{r}\right] !\right) }\left[ p_{n+1-k_{1}}\right] %
\left[ p_{k_{1}-k_{2}}\right] ...\left[ p_{k_{r-1}-k_{r}}\right] \text{,}
\end{equation*}
where the sum is over the set of sequences of integers $k=\left(
k_{1},k_{2},...,k_{r}\right) $, with $r\in \mathbb{N}^{\ast }$, and $n\geq
k_{1}>k_{2}>...>k_{r-1}>k_{r}=0$. This is precisely Eq. $\left( 4.3\right) $
for $n+1$.
\end{proof}

\bigskip

\begin{theorem}
\bigskip For $n\geq 1$, 
\begin{equation}
e_{n}=\sum\limits_{\left| \lambda \right| =n}\varepsilon _{\lambda }\left[
z_{\lambda }\right] _{q}^{-1}\left[ p_{\lambda }\right] _{q}\text{,} 
\tag{4.5}
\end{equation}
where, for any partition $\lambda =\left( \lambda _{1},\lambda
_{2},...,\lambda _{r}\right) $ of $n$: 
\begin{equation}
\left[ z_{\lambda }\right] _{q}=\left( \sum\nolimits_{u}\left( \left[ n%
\right] \left[ n-u_{1}\right] \left[ n-u_{1}-u_{2}\right] ...\left[ u_{r}%
\right] \right) ^{-1}\right) ^{-1}\text{,}  \tag{4.6}
\end{equation}
where the sum extends over the set$\ $of distinct permutations $u$\ of $%
\lambda $.
\end{theorem}

Note that when $q=1$, $\left[ p_{\lambda }\right] =p_{\lambda }$. As the $%
p_{\lambda }$ form a basis of $\Lambda _{K}$, we deduce from the comparison
of $\left( 4.5\right) $ and $\left( 1.3\right) $ that $\left[ z_{\lambda }%
\right] _{q=1}=z_{\lambda }$. This proves that $\left[ z_{\lambda }\right]
_{q}$ is a $q$-analog of $z_{\lambda }$ and that $\left( 4.5\right) $ is a
q-analog of $\left( 1.3\right) $.

\begin{proof}
\begin{proof}
We start from Lemma 4.1. For each sequence $k=\left(
k_{1},k_{2},...,k_{r}\right) $, with $r\in \mathbb{N}^{\ast }$, and $n-1\geq
k_{1}>k_{2}>...>k_{r-1}>k_{r}=0$, we associate the sequence $u=\left(
u_{1},u_{2},...,u_{r}\right) $ defined by:
\end{proof}

$\qquad \qquad \left\{ 
\begin{array}{c}
u_{1}=n-k_{1} \\ 
u_{2}=k_{1}-k_{2} \\ 
... \\ 
u_{r}=k_{r-1}-k_{r}
\end{array}
\right. $

We clearly have for $1\leq i\leq r$,$\ \ u_{i}>0$, and $%
u_{1}+u_{2}+...+u_{r}=n-k_{r}=n$. So, $u$ is a strict composition of $n$.
Moreover, the sequences $k$ in the sum of $\left( 4.3\right) $ and the
strict compositions $u$ of $n$ are in bijection. Indeed, to each strict
composition $u=\left( u_{1},u_{2},...,u_{r}\right) $ of $n$, we can
associate the sequence $k=\left( k_{1},k_{2},...,k_{r}\right) $, defined by

\begin{equation}
\left\{ 
\begin{array}{c}
k_{1}=n-u_{1} \\ 
k_{2}=k_{1}-u_{2}=n-\left( u_{1}+u_{2}\right)  \\ 
... \\ 
k_{r}=k_{r-1}-k_{r}=n-\left( u_{1}+u_{2}+...+u_{r}\right) =0
\end{array}
\right. \text{.}  \tag{4.7}
\end{equation}
This sequence satisfies the inequalities $n-1\geq
k_{1}>k_{2}>...>k_{r-1}>k_{r}=0$. Therefore, $(4.3)$ can be rewritten, by
replacing $k_{i}$ with their expressions in terms of $u_{i}$ and dividing by 
$\left[ n\right] !$, as: 
\begin{equation*}
e_{n}=\sum\nolimits_{u}\left( -1\right) ^{n-r}\left( \left[ n\right] \left[
n-u_{1}\right] \left[ n-\left( u_{1}+u_{2}\right) \right] ...\left[ u_{r}%
\right] \right) ^{-1}\left[ p_{u_{1}}\right] \left[ p_{u_{2}}\right] ...%
\left[ p_{u_{r}}\right] \text{,}
\end{equation*}
where the sum is over the strict compositions $u$ of $n$. For any strict
composition $u=\left( u_{1},u_{2},...,u_{r}\right) $ of $n$, let $\lambda
_{u}$ denote the partition of $n$, composed of the terms of $u$ arranged in
non-increasing order. Given a partition $\lambda =\left( \lambda
_{1},\lambda _{2},...,\lambda _{r}\right) $ of $n$, all the strict
compositions $u$\ of $n$ such that $\lambda _{u}=\lambda $ are obtained in
one and only one way, by taking all the distinct permutations of $\left(
\lambda _{1},\lambda _{2},...,\lambda _{r}\right) $. Moreover if $\lambda
_{u}=\lambda $, then $\left( -1\right) ^{n-r}=\left( -1\right) ^{\left|
\lambda \right| -l(\lambda )}=\varepsilon _{\lambda }$ and $\left[ p_{u_{1}}%
\right] \left[ p_{u_{2}}\right] ...\left[ p_{u_{r}}\right] =\left[
p_{\lambda }\right] $. We therefore have $e_{n}=\sum\nolimits_{\left|
\lambda \right| =n}\varepsilon _{\lambda }\left[ p_{\lambda }\right]
\sum\nolimits_{u}\left( \left[ n\right] \left[ n-u_{1}\right] \left[
n-\left( u_{1}+u_{2}\right) \right] ...\left[ u_{r}\right] \right) ^{-1}$($%
u\,$distinct permutations of $\lambda $), which gives $\left( 4.5\right) $
with $\left( 4.6\right) $.
\end{proof}

\bigskip

\begin{corollary}
\bigskip For $\lambda =\emptyset $ , set\ $\left| \lambda \right| =0$ and $%
l\left( \lambda \right) =0.$ Then 
\begin{equation*}
\mathbf{E}\left( t\right) =\sum\nolimits_{n\geq
0}e_{n}t^{n}=\sum\nolimits_{\lambda }\varepsilon _{\lambda }\left[
z_{\lambda }\right] _{q}^{-1}\left[ p_{\lambda }\right] _{q}t^{\left|
\lambda \right| }\text{,}
\end{equation*}
which is a q-analog of $\left( 4.1\right) $.
\end{corollary}

\begin{proof}
\bigskip immediate with the definition of $E(t)$ and $\left( 4.5\right) $.
\end{proof}

\subsection{\protect\bigskip Case of $h_{n}$}

\bigskip We proceed in the same way as for $e_{n}$.

\begin{lemma}
\bigskip For $n\geq 1$, we have: 
\begin{equation}
\left[ n\right] !h_{n}=\sum\limits_{k}q^{k_{1}+k_{2}+...+k_{r}}\dfrac{\left[
n-1\right] !}{\left[ k_{1}\right] \left[ k_{2}\right] ...\left[ k_{r-1}%
\right] \left( \left[ k_{r}\right] !\right) }\left[ p_{n-k_{1}}\right] \left[
p_{k_{1}-k_{2}}\right] ...\left[ p_{k_{r-1}-k_{r}}\right] \text{,}  \tag{4.8}
\end{equation}
where the sum extends over all sequences of integers $k=\left(
k_{1},k_{2},...,k_{r}\right) $, with $r\in \mathbb{N}^{\ast }$ and $n-1\geq
k_{1}>k_{2}>...>k_{r-1}>k_{r}=0$.
\end{lemma}

\begin{proof}
By induction. This is true for $n=1$. Indeed $\left[ 1\right] !h_{1}=p_{1}$,
and the set of $k$ in the sum of $\left( 4.8\right) $ is reduced to $r=1$
and $k_{1}=0$, so this sum is also $q^{0}\left( 0!/0!\right) \left[ p_{1}%
\right] =p_{1}$. Assuming true $\left( 4.8\right) $ for $n\geq j\geq 1$,
let's prove it for $n+1$. Multiplying the $(n+1)$th equation $\left(
3.7\right) $ by $\left[ n\right] !$ gives 
\begin{equation}
\left[ n+1\right] !h_{n+1}=\sum\limits_{j=0}^{n}\dfrac{\left[ n\right] !}{%
\left[ j\right] !}q^{j}\left[ p_{n+1-j}\right] \left[ j\right] !h_{j} 
\tag{4.9}
\end{equation}
After the same transformation as before and with the induction hypothesis on
the $h_{i}$, $i\leq n$, we obtain$:$%
\begin{equation*}
\left[ n+1\right] !h_{n+1}=\sum\limits_{k}q^{k_{1}+k_{2}+...+k_{r}}\dfrac{%
\left[ n\right] !}{\left[ k_{1}\right] \left[ k_{2}\right] ...\left[ k_{r-1}%
\right] \left( \left[ k_{r}\right] !\right) }\left[ p_{n+1-k_{1}}\right] %
\left[ p_{k_{1}-k_{2}}\right] ...\left[ p_{k_{r-1}-k_{r}}\right] \text{,}
\end{equation*}
where $k=\left( k_{1},k_{2},...,k_{r}\right) $ and $n\geq
k_{1}>k_{2}>...>k_{r}=0$) which is $\left( 4.8\right) $ for $n+1$.
\end{proof}

\begin{theorem}
\bigskip \bigskip For $n\geq 1$, 
\begin{equation}
h_{n}=\sum\limits_{\left| \lambda \right| =n}(q^{\left| \lambda \right|
-l\left| \lambda \right| }\left[ z_{\lambda }\right] _{q^{-1}})^{-1}\left[
p_{\lambda }\right] _{q}\text{,}  \tag{4.10}
\end{equation}
with for any partition $\lambda =\left( \lambda _{1},\lambda
_{2},...,\lambda _{r}\right) $ of $n$: 
\begin{equation}
\left( \left[ z_{\lambda }\right] _{q^{-1}}\right)
^{-1}=\sum\limits_{u}\left( \left[ n\right] _{q^{-1}}\left[ n-u_{1}\right]
_{q^{-1}}\left[ n-u_{1}-u_{2}\right] _{q^{-1}}...\left[ u_{r}\right]
_{q^{-1}}\right) ^{-1}\text{,}  \tag{4.11}
\end{equation}
where the sum extends over the set$\ $of distinct permutations $u$\ of $%
\lambda $.
\end{theorem}

As before, we deduce from the comparison of $\left( 4.10\right) $ and $%
\left( 1.4\right) $ that $\left( q^{\left| \lambda \right| -l(\lambda )}%
\left[ z_{\lambda }\right] _{q^{-1}}\right) _{q=1}=z_{\lambda }$. This
proves that $\left( q^{\left| \lambda \right| -l(\lambda )}\left[ z_{\lambda
}\right] _{q^{-1}}\right) $ is another $q$-analog of $z_{\lambda }$ and that 
$\left( 4.10\right) $ is a q-analog of $\left( 1.4\right) $.

\begin{proof}
\bigskip From Lemma 4.4, we make the same change of index as previously,
which gives the same transfomation of 
\begin{equation*}
\dfrac{1}{\left[ k_{1}\right] \left[ k_{2}\right] ...\left[ k_{r-1}\right]
\left( \left[ k_{r}\right] !\right) }\left[ p_{n-k_{1}}\right] \left[
p_{k_{1}-k_{2}}\right] ...\left[ p_{k_{r-1}-k_{r}}\right] 
\end{equation*}
in terms of the $u_{i}$. Moreover, Eq. $\left( 4.7\right) $ and the fact
that $\left| u\right| =n$ imply: 
\begin{equation*}
k_{1}+k_{2}+...+k_{r}=nr-r\left( u_{1}+u_{2}+...+u_{r}\right)
+\sum\nolimits_{j=1}^{r}\left( j-1\right) u_{j}=n(u)\text{.}
\end{equation*}
By dividing Eq. $\left( 4.8\right) $ by $\left[ n\right] !$, we obtain, in
the same way as for $e_{n}$: 
\begin{equation*}
h_{n}=\sum\nolimits_{\left| \lambda \right| =n}\left[ p_{\lambda }\right]
\sum\nolimits_{u}q^{n\left( u\right) }\left( \left[ n\right] \left[ n-u_{1}%
\right] \left[ n-\left( u_{1}+u_{2}\right) \right] ...\left[ u_{r}\right]
\right) ^{-1}\text{,}
\end{equation*}
where the rightmost sum is over distinct permutations $u$\ of $\lambda $.
Let's now show that this sum is equal to $(q^{\left| \lambda \right|
-l\left( \lambda \right) }\left[ z_{\lambda }\right] _{q^{-1}})^{-1}$. Eq. $%
\left( 4.6\right) $ gives 
\begin{equation*}
(\left[ z_{\lambda }\right] _{q^{-1}})^{-1}=\sum_{u}(\left[ n\right]
_{q^{-1}}\left[ n-u_{1}\right] _{q^{-1}}...\left[ u_{r}\right]
_{q^{-1}})^{-1}\text{.}
\end{equation*}
Since $\left[ m\right] _{q^{-1}}=\left[ m\right] _{q}q^{-m+1}$ for any
integer $m$, we obtain: 
\begin{equation*}
(q^{\left| \lambda \right| -l\left( \lambda \right) }\left[ z_{\lambda }%
\right] _{q^{-1}})^{-1}=\sum\nolimits_{u}q^{\Xi }\left( \left[ n\right] %
\left[ n-u_{1}\right] \left[ n-\left( u_{1}+u_{2}\right) \right] ...\left[
u_{r}\right] \right) ^{-1}\text{,}
\end{equation*}
with 
\begin{equation*}
\Xi =l(\lambda )-\left| \lambda \right|
+n-1+n-u_{1}-1+n-(u_{1}+u_{2)}-1+...+n-(u_{1}+...+u_{r-1})-1=n\left(
u\right) \text{,}
\end{equation*}
as seen before.
\end{proof}

\begin{corollary}
\bigskip \bigskip For $\lambda =\emptyset $, set\ $\left| \lambda \right| =0$
and $l\left( \lambda \right) =0$, then 
\begin{equation*}
H\left( t\right) =\sum\nolimits_{n\geq 0}h_{n}t^{n}=\sum\nolimits_{\lambda
}(q^{\left| \lambda \right| -l(\lambda )}\left[ z_{\lambda }\right]
_{q^{-1}})^{-1}\left[ p_{\lambda }\right] _{q}t^{\left| \lambda \right| }%
\text{,}
\end{equation*}
which is a q-analog of $(4.2)$.
\end{corollary}

\section{\protect\bigskip $E(t)$ as $q$-composition}

Let $P(t)$ be given by $\left( 1.5\right) $. Note that this notation, used
in [2, p.58], differs from the one given by [13, Eq. (2.10) p.23]. The goal
is to find a $q-$analog of the following classical formula (see, for example
[2, Eq. $\left( 3.18\right) $]): 
\begin{equation}
E(t)=\exp (-P(-t))  \tag{5.1}
\end{equation}

To achieve this, we use a formula introduced by Gessel in $\left[ 8\right] $%
. First the definitions introduced in [8] are recalled and slightly
modified, to strengthen the analogy with the classical case. In sections 5
and 6, all the fps belong to $A\left[ \left[ t\right] \right] $, where $A$
is a commutative algebra over a field containing $\mathbb{Q}\left( q\right) $%
.

\subsection{\protect\bigskip q-composition of Gessel}

\begin{definition}
\bigskip (Def. 3.1 from\ [8], modified) Let $F(t)$ be a fps with $F(0)=0$,
and $k\in \mathbb{N}$. The $q$-analog $F^{\left[ k\right] _{q}}$ of $F^{k}$,
is defined as follows: 
\begin{equation}
F^{\left[ 0\right] _{q}}=1\text{ and for }k>0\;\;D_{q}F^{\left[ k\right]
_{q}}=\left[ k\right] F^{\left[ k-1\right] _{q}}D_{q}F\text{.}  \tag{5.2}
\end{equation}
\end{definition}

If there is no ambiguity, the index $q$ will be implicit in writting $F^{%
\left[ k\right] }$. It is easy to verify that:

$\left( i\right) $ For $q=1$ $F^{\left[ k\right] _{q=1}}=F^{k}$. This is the
advantage of this definition compared to the one given in [8].

$\left( ii\right) $ \ $F^{\left[ 1\right] }=F$

$\left( iii\right) $ In general $F^{\left[ j\right] }F^{\left[ k\right]
}\neq F^{\left[ j+k\right] }$

$\left( iii\right) $ $t^{\left[ k\right] }=t^{k}$

\bigskip If we set: 
\begin{equation}
F^{\left[ k\right] }(t)=\sum\limits_{n=0}f_{n,k}\dfrac{t^{n}}{\left[ n\right]
!}\text{,}  \tag{5.3}
\end{equation}
$\left( 5.2\right) $ gives by equating coefficients: 
\begin{equation}
f_{n+1,k}=\left[ k\right] \sum\limits_{j=0}^{n}\QATOPD[ ] {n}{j}%
f_{n-j+1,1}f_{j,k-1}\text{.}  \tag{5.4}
\end{equation}
From this, it is easy to show by induction that $f_{n,k}=0$ for $n<k$.

\begin{definition}
(Modified from [8, Def 3.2]) Let $F(t)$ and $G(t)$ be two fps with $%
G(t)=\sum\nolimits_{n=0}^{\infty }g_{n}\left( t^{n}/\left[ n\right] !\right) 
$ and $F\left( 0\right) =0.$ Then, the $q$-composition $G\left[ F\right] _{q}
$ is defined as: 
\begin{equation}
G\left[ F\right] _{q}=\sum\limits_{k=0}^{\infty }g_{k}\dfrac{F^{\left[ k%
\right] }}{\left[ k\right] !}  \tag{5.5}
\end{equation}
\end{definition}

This definition is valid because $f_{n,k}=0$ for $n<k,$ ensures that the
series in $\left( 5.5\right) $ formally converges in $A\left[ \left[ t\right]
\right] $. Moreover, the two modifications made to the definitions cancel
each other out, resulting in the same composition as in [8]. If $q=1$, $G%
\left[ F\right] _{q=1}$ is the usual composition $G\circ F$. The index $q$
will be omitted if there is no ambiguity in simply writing $G\left[ F\right] 
$.

\begin{proposition}
(Proposition 3.3 of [8]): $D_{q}(G\left[ F\right] )=(D_{q}G)\left[ F\right]
.D_{q}F.$
\end{proposition}

\begin{proof}
$D_{q}(G\left[ F\right] )=\sum\nolimits_{k\geq 0}(g_{k}/\left[ k\right]
!)D_{q}F^{\left[ k\right] }=\sum\nolimits_{k\geq 1}g_{k}\left( D_{q}F\right)
F^{\left[ k-1\right] }/\left[ k-1\right] !=(D_{q}G)\left[ F\right] .D_{q}F$
\end{proof}

Now, let $G$ be the $q$-analog of the exponential function given by (see
[10, Chap.9]):

\begin{equation}
\mathbf{e}_{q}\left( t\right) =\sum\limits_{n\geq 0}\dfrac{t^{n}}{\left[ n%
\right] _{q}!}\text{.}  \tag{5.6}
\end{equation}
Since $D_{q}\mathbf{e}_{q}=\mathbf{e}_{q}$, we have by Proposition 5.3, $%
D_{q}\left( \mathbf{e}_{q}\left[ F\right] \right) =\mathbf{e}_{q}\left[ F%
\right] .D_{q}F$. Equating coefficients give the following recurrence:

\begin{proposition}
\bigskip (Proposition 3.4 of [8]) Let $F\left( t\right)
=\sum\nolimits_{n=1}^{\infty }f_{n}\left( t^{n}/\left[ n\right] !\right) $,
then $\mathbf{e}_{q}\left[ F\right] =\sum\nolimits_{n=0}^{\infty }\gamma
_{n}\left( t^{n}/\left[ n\right] !\right) $ is equivalent to 
\begin{equation}
\gamma _{0}=1\text{ and for }n\geq 0\ \ \gamma _{n+1}=\sum\nolimits_{k=0}^{n}%
\QATOPD[ ] {n}{k}\gamma _{n-k}f_{k+1}  \tag{5.7}
\end{equation}
\end{proposition}

\bigskip

\subsection{\protect\bigskip Application to symmetric functions}

\begin{theorem}
Let the q-analog of $P(t)$ be given by: 
\begin{equation}
P_{q}\left( t\right) =\sum\limits_{n=1}^{\infty }\dfrac{\left[ p_{n}\right] 
}{\left[ n\right] }t^{n}\text{,}  \tag{5.8}
\end{equation}
then, we have the q-analog of $\left( 5.1\right) :$ 
\begin{equation}
E\left( t\right) =\mathbf{e}_{q}\left[ -P_{q}\left( -t\right) \right] _{q}%
\text{.}  \tag{5.9}
\end{equation}
\end{theorem}

\begin{proof}
Let $A$ be the algebra of symmetric functions $\Lambda _{K}$. Equation $%
\left( 4.4\right) $ can be rewritten with $k=n-j$%
\begin{equation*}
\left[ n+1\right] !e_{n+1}=\sum\limits_{k=0}^{n}\left( -1\right) ^{k}\dfrac{%
\left[ n\right] !}{\left[ n-k\right] !\left[ k\right] !}\left[ n-k\right]
!e_{n-k}\left[ k\right] !\left[ p_{k+1}\right] 
\end{equation*}
which is equivalent to $\left( 5.7\right) $, with $\gamma _{n}=\left[ n%
\right] !e_{n}$ and $f_{n}=\left[ n-1\right] !\left( -1\right) ^{n-1}\left[
p_{n}\right] $. Therefore, $F(t)=-P_{q}\left( -t\right) $, and by
Proposition 5.4 $\mathbf{e}_{q}\left[ -P_{q}\left( -t\right) \right]
=\sum\nolimits_{n\geq 0}(\left[ n\right] !e_{n}/\left[ n\right] !)t^{n}=E(t)$
\end{proof}

\section{\protect\bigskip $H(t)$ as $q^{\ast }$-composition}

Still using $\left( 1.5\right) $, we seek a $q-$analog of the classical
formula (see for example [2, Eq. $\left( 3.18\right) $]):

\begin{equation}
H(t)=\exp (P(t))\text{.}  \tag{6.1}
\end{equation}

To achieve this, we introduce another analog of the composition.

\subsection{\protect\bigskip $q^{\ast }$-composition}

\begin{definition}
\bigskip Let $F\left( t\right) $ be a fps with $F(0)=0$ and $k\in \mathbb{N}$%
. The $q^{\ast }$-analog $F^{\left[ k\right] _{q}^{\ast }}$ of $F^{k}$, is
defined by: 
\begin{equation}
F^{\left[ 0\right] _{q}^{\ast }}=1\text{ and for }k>0\;\;D_{q}F^{\left[ k%
\right] _{q}^{\ast }}=\left[ k\right] _{q}q^{-\left( k-1\right) }\left( F^{%
\left[ k-1\right] _{q}^{\ast }}\left( qt\right) \right) D_{q}F\text{.} 
\tag{6.2}
\end{equation}
\end{definition}

\bigskip If there is no ambiguity, the index $q$ will be implied by simply
writting $F^{\left[ k\right] ^{\ast }}$. It is easy to check that:

$\left( i\right) $ For $q=1$ $F^{\left[ k\right] ^{\ast }}=F^{k}$

$\left( ii\right) $ \ $F^{\left[ 1\right] ^{\ast }}=F$

$\left( iii\right) $ in general $F^{\left[ j\right] ^{\ast }}F^{\left[ k%
\right] ^{\ast }}\neq F^{\left[ j+k\right] ^{\ast }}$

$\left( iii\right) $ $t^{\left[ k\right] ^{\ast }}=t^{k}$

\begin{proposition}
\bigskip If we set $F^{\left[ k\right] ^{\ast
}}(t)=\sum\nolimits_{n=0}^{\infty }f_{n,k}t^{n}/\left[ n\right] !$, than 
\begin{equation}
f_{n+1,k}=\left[ k\right] q^{-\left( k-1\right) }\sum\limits_{j=0}^{n}%
\QATOPD[ ] {n}{j}\;f_{n-j+1,1}\;f_{j,k-1}\;q^{j}  \tag{6.3}
\end{equation}
and $f_{n,k}=0$ for $n<k$
\end{proposition}

\begin{proof}
\bigskip On the one hand we have 
\begin{equation*}
D_{q}F^{\left[ k\right] ^{\ast }}=\sum\nolimits_{n=1}^{\infty }f_{n,k}\dfrac{%
t^{n-1}}{\left[ n-1\right] !}=\sum\nolimits_{n=0}^{\infty }f_{n+1,k}\dfrac{%
t^{n}}{\left[ n\right] !}\text{.}
\end{equation*}
On the other hand, $D_{q}F^{\left[ k\right] ^{\ast }}=\left[ k\right]
q^{-\left( k-1\right) }(F^{\left[ k-1\right] ^{\ast }}\left( qt\right)
)D_{q}F=\left[ k\right] q^{-\left( k-1\right) }\left(
\sum\nolimits_{n=1}^{\infty }f_{n,k-1}q^{n}t^{n}/\left[ n\right] !\right)
\left( \sum\nolimits_{n=0}^{\infty }f_{n+1,1}t^{n}/\left[ n\right] !\right) $%
. Equating coefficients gives $\left( 6.3\right) $. Then, by induction on $k$%
, we deduce that $f_{n,k}=0$ for $n<k$, starting from $f_{0,1}=F(0)=0$.
\end{proof}

\begin{definition}
\bigskip Let $G(t)$ and $F(t)$ be two fps with $G(t)=\sum\nolimits_{n=0}^{%
\infty }g_{n}\left( t^{n}/\left[ n\right] !\right) $ and $F\left( 0\right)
=0.$ Then, the $q^{\ast }$-composition $G\left[ F\right] _{q}^{\ast }$ is
defined as 
\begin{equation}
G\left[ F\right] _{q}^{\ast }=\sum\limits_{k=0}^{\infty }g_{k}\dfrac{F^{%
\left[ k\right] ^{\ast }}}{\left[ k\right] _{q}!}\text{.}  \tag{6.4}
\end{equation}
\end{definition}

This definition makes sense because $f_{n,k}=0$ for $n<k,$ so the series in $%
\left( 6.4\right) $ converges formally in $A\left[ \left[ t\right] \right] $%
. If $q=1$, $G\left[ F\right] _{q=1}^{\ast }$ is the usual composition $%
G\circ F$. The index $q$ will be omitted if there is no ambiguity, simply
writting $G\left[ F\right] ^{\ast }$.

\begin{proposition}
$\ D_{q}(G\left[ F\right] ^{\ast })(t)=((D_{q}G)\left[ q^{-1}F\right] ^{\ast
})\left( qt\right) .(D_{q}F)$.
\end{proposition}

\bigskip We find the usual chain rule if $q=1$. For the proof, we use the
following lemma.

\begin{lemma}
If $F_{\theta }\left( t\right) =\theta \left( q\right) F\left( t\right) $
where $\theta \left( q\right) $ $\in \mathbb{Q}\left( q\right) $ (i.e., $%
\theta \left( q\right) $ independant of $t$), then, 
\begin{equation*}
F_{\theta }^{\left[ k\right] ^{\ast }}=\left( \theta \left( q\right) \right)
^{k}F^{\left[ k\right] ^{\ast }}\text{.}
\end{equation*}
\end{lemma}

\begin{proof}
By induction. \bigskip This is clearly true for $k=0$, both sides equaling
1. Assume the equation of the lemma is true for $k-1$, then 
\begin{eqnarray*}
D_{q}\left( F_{\theta }^{\left[ k\right] ^{\ast }}\right)  &=&\left[ k\right]
q^{-\left( k-1\right) }\left( F_{\theta }^{\left[ k-1\right] ^{\ast }}\left(
qt\right) \right) D_{q}F_{\theta }=\left[ k\right] q^{-\left( k-1\right)
}\left( \theta \left( q\right) \right) ^{k-1}\left( F^{\left[ k-1\right]
^{\ast }}\left( qt\right) \right) \theta \left( q\right) D_{q}F \\
&=&\left( \theta \left( q\right) \right) ^{k}D_{q}F^{\left[ k\right] ^{\ast
}}=D_{q}\left( \left( \theta \left( q\right) \right) ^{k}F^{\left[ k\right]
^{\ast }}\right) 
\end{eqnarray*}
and moreover the two fps $F_{\theta }^{\left[ k\right] ^{\ast }}$ et $\left(
\theta \left( q\right) \right) ^{k}F^{\left[ k\right] ^{\ast }}$ are $0$ for 
$t=0$.
\end{proof}

\bigskip

\begin{proof}[Proof of Proposition 6.4]
\bigskip 
\begin{equation*}
D_{q}(G\left[ F\right] ^{\ast })=\sum\nolimits_{k\geq 0}\dfrac{g_{k}}{\left[
k\right] !}D_{q}F^{\left[ k\right] ^{\ast }}=\sum\nolimits_{k\geq 1}\dfrac{%
g_{k}}{\left[ k\right] !}\left[ k\right] q^{-\left( k-1\right) }\left( F^{%
\left[ k-1\right] ^{\ast }}\left( qt\right) \right) D_{q}F\text{,}
\end{equation*}
hence, by Lemma 6.5 
\begin{equation*}
D_{q}(G\left[ F\right] ^{\ast })=\left( \sum\nolimits_{k\geq 1}\dfrac{g_{k}}{%
\left[ k-1\right] !}\left( q^{-1}F\right) ^{\left[ k-1\right] ^{\ast
}}\right) \left( qt\right) .D_{q}F=\left( (D_{q}G)\left[ q^{-1}F\right]
^{\ast }\right) \left( qt\right) .(D_{q}F)\text{.}
\end{equation*}
\end{proof}

\bigskip

Let $G$ now represent the second q-analog of the exponential function given
by (see [10, Chap. 9]):

\begin{equation}
\mathbf{E}_{q}\left( t\right) =\sum\limits_{n\geq 0}q^{\binom{n}{2}}\dfrac{%
t^{n}}{\left[ n\right] _{q}!}\text{.}  \tag{6.5}
\end{equation}

\begin{proposition}
\bigskip Let $F\left( t\right) =\sum\nolimits_{n=1}^{\infty }f_{n}\left(
t^{n}/\left[ n\right] !\right) $ then $\mathbf{E}_{q}\left[ F\right] ^{\ast
}=\sum\nolimits_{n=0}^{\infty }\gamma _{n}^{\ast }\left( t^{n}/\left[ n%
\right] !\right) $ is equivalent to 
\begin{equation}
\text{ }\gamma _{0}^{\ast }=1\text{ \ and for }n\geq 0\text{\ \ }\gamma
_{n+1}^{\ast }=\sum\limits_{k=0}^{n}\QATOPD[ ] {n}{k}q^{n-k}\gamma
_{n-k}^{\ast }f_{k+1}\text{.}  \tag{6.6}
\end{equation}
\end{proposition}

\begin{proof}
\bigskip $\mathbf{E}_{q}$ satisfies$\ D_{q}\mathbf{E}_{q}(t)=\mathbf{E}%
_{q}\left( qt\right) =\sum\nolimits_{n\geq 0}q^{\binom{n}{2}+n}t^{n}/\left[ n%
\right] !$. Therefore, by Proposition \ 6.4: \ \ 
\begin{equation}
D_{q}\left( \mathbf{E}_{q}\left[ F\right] ^{\ast }\right) =(D_{q}\mathbf{E}%
_{q})\left[ q^{-1}F\right] ^{\ast }\left( qt\right) .D_{q}F\text{.} 
\tag{6.7}
\end{equation}

We have on the left side of $\left( 6.7\right) $,\ $\sum\nolimits_{n\geq
0}\gamma _{n+1}^{\ast }t^{n}/\left[ n\right] !$, and on the right side of $%
\left( 6.7\right) $, with Lemma 6.5: 
\begin{equation*}
\sum\nolimits_{n\geq 0}\dfrac{q^{\binom{n}{2}+n}}{\left[ n\right] !}%
q^{-n}\left( F\right) ^{\left[ n\right] ^{\ast }}\left( qt\right) .D_{q}F=%
\mathbf{E}_{q}\left[ F\right] ^{\ast }\left( qt\right)
.D_{q}F=\sum\nolimits_{k\geq 0}\gamma _{k}^{\ast }q^{k}\dfrac{t^{k}}{\left[ k%
\right] !}\sum\nolimits_{l\geq 0}f_{l+1}\dfrac{t^{l}}{\left[ l\right] !}%
\text{.}
\end{equation*}
Equating coefficients in $\left( 6.7\right) $ gives Recurrence $\left(
6.6\right) .$ Conversely, these equalities between coefficients give $%
\mathbf{E}_{q}\left[ F\right] ^{\ast }=\sum\nolimits_{n\geq 0}\gamma
_{n}^{\ast }t^{n}/\left[ n\right] !$.
\end{proof}

\bigskip

\textbf{Non-inversion in permutation}

In $\left[ 8\right] $, Gessel used Proposition 5.4 for the enumeration of
inversions of permutations. We will succintly show the use of Proposition $%
6.6$ for the enumeration of non-inversions of permutations. First, we
summarize the notations of Gessel, which we extend to non-inversions,
referring to $\left[ 8\right] $ for more details.

Let's denote by $a_{1}a_{2}...a_{n}$ a permutation of a set $E\subset 
\mathbb{N}$ with $\left| E\right| =n$. A permutation is basic if it begins
with its greatest element. We denote by $S_{n}$ and $B_{n}$ the sets of
permutations and basic permutations on \textbf{n}, respectively. The content
of the permutation $\pi =a_{1}a_{2}...a_{n}$ is $con(\pi )=\left\{
a_{1},a_{2,}...,a_{n}\right\} $. It is shown in $\left[ 8\right] $ that
there is a bijection between $S_{n}$ and the set of disjoint basic
permutations $\left\{ \beta _{1},\beta _{2},...\beta _{k}\right\} $ whose
union of contents is \textbf{n}. For $\pi \in S_{n}$ the corespondant set $%
\left\{ \beta _{1},\beta _{2},...\beta _{k}\right\} $ is called the basic
decomposition of $\pi $. We call a permutation reduced if it is in $S_{n}$
for some $n\in \mathbb{N}$. To any permutation $\pi =a_{1}a_{2}...a_{n}$ we
may associate a reduced permutation $red\left( \pi \right) $, by replacing
in $\pi $, for each $i=1,2,...,n$, the $i$th smallest element of $\left\{
a_{1},a_{2,}...,a_{n}\right\} $ by $i$.

A function $\omega $ defined on permutations, with values in some comutative
algebra over $\mathbb{Q}$, is multiplicative if for all permutations $\pi $:

$\left( i\right) $ $\omega \left( \pi \right) =\omega \left( red\left( \pi
\right) \right) $

$\left( ii\right) $ If $\left\{ \beta _{1},\beta _{2},...\beta _{k}\right\} $
is the basic decomposition of $\pi $, then

\begin{equation*}
\omega \left( \pi \right) =\omega \left( \beta _{1}\right) \omega \left(
\beta _{2}\right) ...\omega \left( \beta _{k}\right) \text{.}
\end{equation*}

If $V$ is a subset of $n$, we denote by $I_{n}(V)$ and $\overline{I_{n}}(V)$%
, the number of pairs $\left( v,w\right) $ with $v\in V$, $w\in \mathbf{n}/V$%
, $v>w$ and $v<w$, respectively.

For the permutation $\pi =a_{1}a_{2}...a_{n}$, an inversion and a
non-inversion, is a pair $\left( i,j\right) $ with $i<j$, $a_{i}>a_{j}$ and $%
a_{i}<a_{j}$, respectively. We write $I\left( \pi \right) $ and $\overline{I}%
\left( \pi \right) $ for the number of inversions and non-inversions of $\pi 
$, respectively. The result of Gessel is stated in the following theorem.

\textit{Theorem 5.2 of [8]: Let }$\omega $\textit{\ be a multiplicative
function on permutations. Let }$\gamma _{n}=\sum\nolimits_{\pi \in
S_{n}}\omega \left( \pi \right) q^{I\left( \pi \right) }$\textit{, }$%
f_{n}=\sum\nolimits_{\beta \in B_{n}}\omega \left( \beta \right) q^{I\left(
\beta \right) }$\textit{\ and }$F(t)=\sum\nolimits_{n\geq 1}f_{n}\dfrac{t^{n}%
}{\left[ n\right] _{q}!}$\textit{. Then } 
\begin{equation}
\sum\limits_{n\geq 0}\gamma _{n}\dfrac{t^{n}}{\left[ n\right] _{q}!}=\mathbf{%
e}_{q}\left[ F(t)\right] _{q}  \tag{6.8}
\end{equation}
The proof is obtained combinatorially by proving the relation (5.7) using
the decomposition of permutation into basic permutations. This justifies
referring to it as the $q$-analog of the classic exponential formula (for
this latter, see [15, Chap.5]). Here is the counterpart of this theorem for
non-inversions.

\begin{theorem}
\bigskip Let $\omega $ be a multiplicative function on permutations. Let $%
\overline{\gamma _{n}}=\sum\nolimits_{\pi \in S_{n}}\omega \left( \pi
\right) q^{\overline{I}\left( \pi \right) }$, $\overline{f_{n}}%
=\sum\nolimits_{\beta \in B_{n}}\omega \left( \beta \right) q^{\overline{I}%
\left( \beta \right) }$ and $\overline{F}(t)=\sum\nolimits_{n\geq 1}%
\overline{f_{n}}t^{n}/\left[ n\right] !$ then 
\begin{equation}
\sum\limits_{n\geq 0}\overline{\gamma _{n}}\dfrac{t^{n}}{\left[ n\right]
_{q}!}=\mathbf{E}_{q}\left[ \overline{F}(t)\right] _{q}^{\ast }\text{.} 
\tag{6.9}
\end{equation}
\end{theorem}

\bigskip Let us first show the counterpart of Lemma 5.1 in $\left[ 8\right] $%
.

\begin{lemma}
Let $\overline{Q}\left( n,k\right) =\sum\nolimits_{V}q^{\overline{I_{n}}%
\left( V\right) }$\ where the sum is over all $V\subseteq \mathbf{n}$ with $%
\left| V\right| =k$. Then $\overline{Q}\left( n,k\right) =\QATOPD[ ] {n}{k}%
_{q}$.
\end{lemma}

\begin{proof}
\bigskip Note that $\overline{I_{n}}\left( V\right) =I_{n}\left( \overline{V}%
\right) $, where $\overline{V}$ is the complement of $V$ in $\mathbf{n}.$ So 
$\overline{Q}\left( n,k\right) =\sum\limits_{\overline{V}}q^{I_{n}\left( 
\overline{V}\right) }$, where the sum is over all $\overline{V}\subseteq 
\mathbf{n}$ with $\left| \overline{V}\right| =n-k$. It is therefore $\QATOPD[
] {n}{k}$, according to Lemma 5.1 of [8].\bigskip
\end{proof}

\bigskip

\begin{proof}
We want to prove the recurrence $\left( 6.6\right) $, i.e. 
\begin{equation}
\overline{\gamma _{0}}=1\text{ \ and for }n\geq 0\text{\ \ }\overline{\gamma
_{n+1}}=\sum\limits_{k=0}^{n}\QATOPD[ ] {n}{k}q^{n-k}\overline{\gamma _{n-k}}%
\overline{f_{k+1}}\text{.}  \tag{6.10}
\end{equation}

We give two proofs of this recurrence.

\textbf{a) Combinatorial proof.} It follows the proof of Theorem 5.2 of [8].
Let show that $\QATOPD[ ] {n}{k}q^{n-k}\overline{\gamma _{n-k}}\overline{%
f_{k+1}}$ counts those permutations counted by $\overline{\gamma _{n+1}}$
whose last basic component has length $k+1$. Such a permutation may be
factorised as $\pi =\sigma \beta $ where $\sigma $ is of length $n-k$, and $%
\beta $ is of length $k+1$. The condition that $\beta $ is the last basic
component of $\pi $, is equivalent to the condition that $\beta $ is basic
and $con\left( \beta \right) $ contains $n+1$. Thus to determine $\pi $ we
choose $V=con(\sigma )$ as an arbitrary $\left( n-k\right) $-subset of $%
\mathbf{n}$, and choose $red\left( \sigma \right) \in S_{n-k}$ and $%
red\left( \beta \right) \in B_{k+1}$. It is easily seen that $\overline{I}%
\left( \pi \right) =\overline{I}\left( \sigma \right) +\overline{I}\left(
\beta \right) +\overline{I_{n}}\left( V\right) +n-k$, the term $n-k$ coming
from the pairs $\left( v,n+1\right) $, $v\in V$. The end of the proof is
similar to that of Theorem 5.2 of [8] using Lemma 6.8\ and $\QATOPD[ ] {n}{%
n-k}=\QATOPD[ ] {n}{k}$.

b) \textbf{Semi-combinatorial proof. }It directly uses Theorem 5.2 of [8].
For $\pi \in $ $S_{n}$\ it is clear that $\overline{I}\left( \pi \right) =%
\binom{n}{2}-I\left( \pi \right) $. Therefore $\overline{\gamma _{n}}%
=\sum\nolimits_{\pi \in S_{n}}\omega \left( \pi \right) q^{\binom{n}{2}%
-I\left( \pi \right) }=q^{\binom{n}{2}}\gamma _{n}\left( 1/q\right) $.
Similarly $\overline{f_{n}}=q^{\binom{n}{2}}f_{n}\left( 1/q\right) $. But
according to Recurrence $\left( 5.7\right) $, we have, by replacing $q$ by $%
1/q$ 
\begin{equation*}
\gamma _{n+1}\left( 1/q\right) =\sum\limits_{k=0}^{n}\QATOPD[ ] {n}{k}%
_{1/q}\gamma _{n-k}\left( 1/q\right) f_{k+1}\left( 1/q\right) \text{.}
\end{equation*}

Multiplying by $q^{\binom{n}{2}}$ we obtain $\left( 6.10\right) $, taking
into account that 
\begin{equation*}
\QATOPD[ ] {n}{k}_{1/q}=\QATOPD[ ] {n}{k}_{q}q^{-\left( n-k\right) k}\text{.}
\end{equation*}
\end{proof}

The combinatorial proof given by a) justifies speaking about a second $q$%
-analog of the classic exponential formula.

\subsection{Application to symmetric functions}

\begin{theorem}
\bigskip Let $P_{q}(t)$ given by $\left( 5.8\right) $, then we have the
q-analog of $\left( 6.1\right) $ 
\begin{equation}
H\left( t\right) =\mathbf{E}_{q}\left[ P_{q}\left( t\right) \right]
_{q}^{\ast }\text{.}  \tag{6.11}
\end{equation}
\end{theorem}

\begin{proof}
\bigskip Let $A$ be the algebra of symmetric functions $\Lambda _{K}$ (see
Section 2). Eq. $\left( 4.9\right) $ can be rewritten with $k=n-j$%
\begin{equation*}
\left[ n+1\right] !h_{n+1}=\sum\limits_{k=0}^{n}\dfrac{\left[ n\right] !}{%
\left[ n-k\right] !\left[ k\right] !}q^{n-k}\left[ n-k\right] !h_{n-k}\left[
k\right] !\left[ p_{k+1}\right] \text{.}
\end{equation*}
This is none other than $\left( 6.6\right) $, with $\gamma _{n}^{\ast }=%
\left[ n\right] !h_{n}$ and $f_{n}=\left[ n-1\right] !\left[ p_{n}\right] $.
Therefore, $F(t)=P_{q}\left( t\right) $ and by Proposition 6.6, $\mathbf{E}%
_{q}\left[ P_{q}\right] ^{\ast }\left( t\right) =\sum\nolimits_{n=0}^{\infty
}(\left[ n\right] !h_{n}/\left[ n\right] !)t^{n}=H(t)$
\end{proof}

\textbf{Open problem. }It is possible to prove combinatorially the classic
formula $\left( 6.1\right) $ (see [15, Chap.5, \S 5.2.11])$.$ It would be
interesting to also prove combinatorially its q-analog $\left( 6.11\right) $.

\section{Link between the two $q$-exponential formulas}

\QTP{Body Math}
In the proofs of [4, Lemma 6.1] and Theorem 4.6 of this article, the
substitution of $q$ by $q^{r}$and of $q$ by $q^{-1}$ was performed to obtain
[4, Eq.(6.1)] and Eq. (4.11), respectively. It is worthwile to elaborate on
the justification of these substitutions by revisiting Definition 2.1 of [4]
in the slightly more general context, where $K$ is the field of Laurent
formal series $\mathbb{Q}\left( \left( q\right) \right) $, and $%
F(t)=\sum\nolimits_{n\geq 0}a_{n}t^{n}\in A\left[ \left[ t\right] \right] $,
with $A$ a commutative algebra over $K$. For the sake of clarity and
sufficiency for our needs, we limit ourselves to the base change $\psi
(q)=q^{m}$ where $m\in \mathbb{Z}^{\ast }$, although the justification can
be generalized without difficulty to $\psi \left( q\right) \in \mathbb{Q}%
\left( \left( q\right) \right) $ with $\psi \left( 1\right) =1$. We define
as in [4, Definition 2.1] for $n\in \mathbb{N}$:

\QTP{Body Math}
\begin{equation*}
\left[ n\right] _{\psi }=1+\psi ^{2}+...+\psi ^{n-1}\text{ if }n\neq 0\text{
and }\left[ 0\right] _{\psi }=0\text{,}
\end{equation*}
then 
\begin{equation*}
\left[ n\right] _{\psi }!=\left[ 1\right] _{\psi }\left[ 2\right] _{\psi }...%
\left[ n\right] _{\psi }\text{ , \ \ }\QATOPD[ ] {n}{k}_{\psi }=\dfrac{\left[
n\right] _{\psi }!}{\left[ k\right] _{\psi }!\left[ n-k\right] _{\psi }!}%
\text{,}
\end{equation*}

\QTP{Body Math}
and 
\begin{equation*}
D_{\psi }F(t)=\dfrac{(F(\psi (q)t)-F(t))}{\left( \psi \left( q\right)
-1\right) t}\text{ and for }r\geq 1\text{ }D_{\psi }^{r}=D_{\psi }(D_{\psi
}^{r-1})\text{ with }D_{\psi }^{0}F=F\text{.}
\end{equation*}
It can be verified that $D_{\psi }^{r}t^{n}=\left[ r\right] _{\psi }!\QATOPD[
] {n}{r}_{\psi }$, so, by linearity 
\begin{equation*}
D_{\psi }^{r}F(t)=\left[ r\right] _{\psi }!\sum\nolimits_{n\geq 0}\QATOPD[ 
] {n}{r}_{\psi }a_{n+1}t^{n-r}\text{.}
\end{equation*}

\QTP{Body Math}
If $A=\Lambda _{K},$ we generalize the definition of $\left[ p_{n}^{\left(
r\right) }\right] $ by

\QTP{Body Math}
\begin{equation*}
\dfrac{D_{\psi }^{(r)}E(t)}{E(t)}=\left[ r\right] _{\psi
}!\sum\limits_{n\geq r}\left[ p_{n}^{(r)}\right] _{\psi }\left( -t\right)
^{n-r}\text{.}
\end{equation*}

\QTP{Body Math}
Since all algebraic calculations leading to Proposition 3.1 and to the
results of Sections 3 and 4, involve the variable $q$, only through $\left[ n%
\right] $ or one of its derivatives $\left[ n\right] !$, $\QATOPD[ ] {n}{k}$
and similarly for the definition of $\left[ p_{n}^{\left( r\right) }\right] $%
, these calculations transpose, mutatis mutandis, by replacing $q$ with $%
\psi \left( q\right) $. As an exemple, the transpositions of $\left(
4.10\right) $ and $\left( 4.11\right) $ obtained in this way are: 
\begin{equation*}
h_{n}=\sum\limits_{\left| \lambda \right| =n}(q^{\left| \lambda \right|
-l(\lambda )}\left[ z_{\lambda }\right] _{\psi ^{-1}})^{-1}\left[ p_{\lambda
}\right] _{\psi }\text{ with }\left( \left[ z_{\lambda }\right] _{\psi
}\right) ^{-1}=\sum\limits_{u}\left( \left[ n\right] _{\psi }\left[ n-u_{1}%
\right] _{\psi }\left[ n-u_{1}-u_{2}\right] _{\psi }...\left[ u_{r}\right]
_{\psi }\right) ^{-1}\text{,}
\end{equation*}
where the summation conditions on $u$\ are the same as in ($4.11)$.

Similarly, we can define the powers and the compositions studied in Sections
5 and 6 with $\psi $ instead of $q$. If we define the fps $\mathbf{e}_{\psi
}(t)=\sum\nolimits_{n\geq 0}t^{n}/\left[ n\right] _{\psi }!$ and $\mathbf{E}%
_{\psi }(t)=\sum\nolimits_{n\geq 0}\psi ^{\binom{n}{2}}t^{n}/\left[ n\right]
_{\psi }!$, it is easy to verify that $D_{\psi }\mathbf{e}_{\psi }=\mathbf{e}%
_{\psi }$ and $D_{\psi }\mathbf{E}_{\psi }=\mathbf{E}_{\psi }\left( \psi
t\right) $. Therefore, all the calculations in Sections 5 and 6, including
Propositions 5.4 and 6.6, Theorems 5.5 and 6.7, are valid by replacing $q$
with $\psi $. For example we have $E(t)=\mathbf{e}_{\psi }\left[ -P_{\psi
}\left( -t\right) \right] _{\psi }$ for all $\psi (q)=q^{m}$ with $m\in 
\mathbb{Z}^{\ast }$, where $P_{\psi }(t)=\sum\nolimits_{n\geq 1}\left[ p_{n}%
\right] _{\psi }t^{n}/\left[ n\right] _{\psi }$.

The link between the two compositions of Sections 5 and 6 is now given by
the following theorem.

\begin{theorem}
\bigskip Let $F\left( t\right) $ be a fps with coefficients in $A$, a
commutative algebra over $\mathbb{Q}\left( \left( q\right) \right) $, with $%
F(0)=0$, and let $\psi (q)=q^{m}$ with $m\in \mathbb{Z}^{\ast }$. Then:

i) For all $k\geq 0$, $\ F^{\left[ k\right] _{\psi }^{\ast }}=F^{\left[ k%
\right] _{\psi ^{-1}}}$

ii) If $G(t)$ is a fps wih coefficients in $A$, then\ \ $G\left[ F\right]
_{\psi }^{\ast }=G\left[ F\right] _{\psi ^{-1}}$

iii) $\mathbf{E}_{\psi }\left[ F\right] _{\psi }^{\ast }=\mathbf{e}_{\psi
^{-1}}\left[ F\right] _{\psi ^{-1}}$. In particular, for $m=1$, there is a
connection between the two $q$-exponential formulas, given by: 
\begin{equation}
\mathbf{E}_{q}\left[ F\right] _{q}^{\ast }=\mathbf{e}_{q^{-1}}\left[ F\right]
_{q^{-1}}  \tag{7.1}
\end{equation}
\end{theorem}

\begin{proof}
\bigskip $i)$ Let's define $F=\sum\nolimits_{n\geq 1}f_{n}t^{n}/\left[ n%
\right] _{\psi }!\ $and, for all $k\geq 0$,$\ \ F^{\left[ k\right] _{\psi
}^{\ast }}=\sum\nolimits_{n\geq 1}f_{n,k}^{\ast }t^{n}/\left[ n\right]
_{\psi }!$\ and\ 
\begin{equation}
F^{\left[ k\right] _{\psi ^{-1}}}=\sum\nolimits_{n\geq 1}f_{n,k}^{\prime
}t^{n}/\left[ n\right] _{\psi ^{-1}}!=\sum\nolimits_{n\geq 1}f_{n,k}t^{n}/%
\left[ n\right] _{\psi }!.  \tag{7.2}
\end{equation}

We reason by induction on $k$. By definition, $F^{\left[ 0\right] _{\psi
^{-1}}}=F^{\left[ 0\right] }=1$, and $F^{\left[ 1\right] _{\psi }^{\ast
}}=F^{\left[ 1\right] _{\psi ^{-1}}}=F$, so $f_{n,1}^{\ast }=f_{n,1}=f_{n}$.
By equating the coefficients in $\left( 7.2\right) $ we have $%
f_{n,k}^{\prime }/\left[ n\right] _{\psi ^{-1}}!=\psi ^{\binom{n}{2}%
}f_{n,k}^{\prime }/\left[ n\right] _{\psi }!=f_{n,k}/\left[ n\right] _{\psi
}!$, so $f_{n,k}=\psi ^{\binom{n}{2}}f_{n,k}^{\prime }$ and in particular, $%
f_{n,1}^{\prime }=\psi ^{-\binom{n}{2}}f_{n}$.

Assume that $F^{\left[ k-1\right] _{\psi }^{\ast }}=F^{\left[ k\right]
_{\psi ^{-1}}}$, therefore $f_{n,k-1}^{\ast }=f_{n,k-1}=\psi ^{\binom{n}{2}%
}f_{n,k-1}^{\prime }$ for each $n\geq 0$.

On one hand, 
\begin{equation}
f_{n+1,k}^{\ast }=\left[ k\right] _{\psi }\psi
^{-(k-1)}\sum\nolimits_{j=0}^{n}\QATOPD[ ] {n}{j}_{\psi }f_{n-j+1,1}^{\ast
}f_{j,k-1}^{\ast }\psi ^{j}=\left[ k\right] _{\psi
^{-1}}\sum\nolimits_{j=0}^{n}\QATOPD[ ] {n}{j}_{\psi }f_{n-j+1}f_{j,k-1}\psi
^{j}.  \tag{7.3}
\end{equation}
\qquad On the other hand, \ $f_{n+1,k}^{\prime }=\left[ k\right] _{\psi
^{-1}}\sum\nolimits_{j=0}^{n}\QATOPD[ ] {n}{j}_{\psi
^{-1}}f_{n-j+1,1}^{\prime }f_{j,k-1}^{\prime }$. Then, 
\begin{equation}
f_{n+1,k}=\psi ^{\binom{n+1}{2}}\left[ k\right] _{\psi
^{-1}}\sum\nolimits_{j=0}^{n}\QATOPD[ ] {n}{j}_{\psi }\psi ^{-j(n-j)}\psi ^{-%
\binom{n-j+1}{2}}f_{n-j+1}\psi ^{-\binom{j}{2}}f_{j,k-1}.  \tag{7.4}
\end{equation}

We can check that 
\begin{equation*}
\binom{n+1}{2}-j(n-j)-\binom{n-j+1}{2}-\binom{j}{2}=j\text{,}
\end{equation*}
so $\left( 7.3\right) $ and $\left( 7.4\right) $ give $f_{n+1,k}^{\ast
}=f_{n+1,k}$ and finally $F^{\left[ k\right] _{\psi }^{\ast }}=F^{\left[ k%
\right] _{\psi ^{-1}}}$.

$ii)$ Let's define $G(t)=\sum\nolimits_{k\geq 0}g_{k}/\left[ k\right] _{\psi
}!=\sum\nolimits_{k\geq 0}g_{k}^{\prime }t^{k}/\left[ k\right] _{\psi
^{-1}}! $, so as above, $g_{k}=\psi ^{\binom{k}{2}}g_{k}^{\prime }$. Then, 
\begin{equation*}
G\left[ F\right] _{\psi }^{\ast }=\sum\nolimits_{k\geq 0}g_{k}\dfrac{F^{%
\left[ k\right] _{\psi }^{\ast }}}{\left[ k\right] _{\psi }!}%
=\sum\nolimits_{k\geq 0}g_{k}^{\prime }\psi ^{\binom{k}{2}}\dfrac{F^{\left[ k%
\right] _{\psi ^{-1}}}}{\left[ k\right] _{\psi }!}=\sum\nolimits_{k\geq
0}g_{k}^{\prime }\dfrac{F^{\left[ k\right] _{\psi ^{-1}}}}{\left[ k\right]
_{\psi ^{-1}}!}=G\left[ F\right] _{\psi ^{-1}}\text{.}
\end{equation*}

$iii)$ It is easy to verify that \ $\mathbf{E}_{\psi }\left( t\right) =%
\mathbf{e}_{\psi ^{-1}}(t)$. Therefore, by replacing in $ii)$: 
\begin{equation*}
\mathbf{E}_{\psi }\left[ F\right] _{\psi }^{\ast }=\mathbf{E}_{\psi }\left[ F%
\right] _{\psi ^{-1}}=\mathbf{e}_{\psi ^{-1}}\left[ F\right] _{\psi ^{-1}}.
\end{equation*}
The particular case is deduced with $\psi \left( q\right) =q$.
\end{proof}

\section{\protect\bigskip infinite q-product for $E(t)$ and $H(t)$}

In this section, formal power series belong to $B\left[ \left[ q,t\right] %
\right] $ with $B=\mathbb{Q}\left[ \left[ Y\right] \right] $, where $Y$ is a
sequence of indeterminates (distinct from $q$ and $t$). If $F\in B\left[ %
\left[ q,t\right] \right] =B\left[ \left[ q\right] \right] \left[ \left[ t%
\right] \right] $, it is necessary to distinguish the order of $F$ as an fps
in $\left( q,t\right) $,\ denoted $\omega \left( F\right) $, from that of $F$%
\ as a fps in $t$ (and coefficients in $B\left[ \left[ q\right] \right] $),
denoted $\omega _{t}\left( F\right) $.

The following definition, whose equivalence of forms arises from the
previous sections, will be useful.

\begin{definition}
\bigskip The fps $p_{q}\left( t\right) \in \Lambda _{K}\left[ \left[ t\right]
\right] $ is defined in four equivalent ways:

i) $p_{q}\left( t\right) =D_{q}P_{q}\left( t\right) $ where $P_{q}$ is given
by $\left( 5.8\right) $

ii) $p_{q}(t)=\sum\nolimits_{n=0}^{\infty }\left[ p_{n+1}\right] t^{n}$

iii) $p_{q}\left( -t\right) =D_{q}E(t)/E(t)$ where $E(t)=\sum\nolimits_{n%
\geq 0}e_{n}t^{n}$

iiii) $p_{q}(t)=D_{q}H(t)/H(qt)$ where $H(t)=\sum\nolimits_{n\geq
0}h_{n}t^{n}$
\end{definition}

\subsection{\protect\bigskip Case of $E(t)$}

\bigskip The tool for this case is still a proposition from [8], which we
generalize in the following way

\begin{proposition}
\bigskip (from Proposition 3.5 of [8]) Let $F\in B\left[ \left[ q,t\right] %
\right] $ with $\omega _{t}\left( F\right) \geq 1$,$\ $and let $\psi \left(
q\right) =q^{m}$ with $m\in \mathbb{N}^{\ast }$. Then::

$i)$ The fps $G=\mathbf{e}_{\psi }\left[ F\right] _{\psi }$ belongs to $B%
\left[ \left[ q,t\right] \right] $ and $\omega _{t}\left( G-1\right) \geq 1$

$ii)$ We have 
\begin{equation}
\mathbf{e}_{\psi }\left[ F\right] _{\psi }=\prod\limits_{k=0}^{\infty
}(1-\left( 1-\psi \right) \psi ^{k}t\left( D_{\psi }F\right) \left( \psi
^{k}t\right) )^{-1}\text{,}  \tag{8.1}
\end{equation}
where the convergence of the product is that of formal convergence in $B%
\left[ \left[ q,t\right] \right] $.
\end{proposition}

\begin{proof}
$i)$ Let $F(t)=\sum\nolimits_{n\geq 1}f_{n}t^{n}/\left[ n\right] _{\psi }!$.
Clearly $B\left[ \left[ q\right] \right] \subset \mathbb{Q}\left( \left(
q\right) \right) \left[ \left[ Y\right] \right] $ so according to
Proposition 5.4 and Section 7 with $A=\mathbb{Q}\left( \left( q\right)
\right) \left[ \left[ Y\right] \right] $, the fps $\mathbf{e}_{\psi }\left[ F%
\right] _{\psi }=\sum\nolimits_{n=0}^{\infty }\gamma _{n}\left( t^{n}/\left[
n\right] _{\psi }!\right) $ is well-defined in $A\left[ \left[ t\right] %
\right] $, and by $\left( 5.7\right) $, $\gamma _{0}=1$ and for $n\geq 0$,$\
\ \gamma _{n+1}=\sum\nolimits_{k=0}^{n}\QATOPD[ ] {n}{k}_{\psi }\gamma
_{n-k}f_{k+1}$. Due to the assumption on $F$ and $\psi $, we have $%
f_{n}/\left( \left[ n\right] _{\psi }!\right) $ $\in B\left[ \left[ q\right] %
\right] $ and $\left[ n\right] _{\psi }!\in \mathbb{Q}\left[ \left[ q\right] %
\right] $, so $f_{n}\in B\left[ \left[ q\right] \right] $. Furthermore, $%
\QATOPD[ ] {n}{k}_{\psi }$is a polynomial in $\psi $, so we deduce by
induction that the coefficients $\gamma _{n}$ and also $\gamma _{n}/\left[ n%
\right] _{\psi }!$ belong to $B\left[ \left[ q\right] \right] $. Thus, $%
G\left( t\right) =\mathbf{e}_{\psi }\left[ F\right] _{\psi }\in B\left[ %
\left[ q,t\right] \right] $, and since $\gamma _{0}=1$, we have $\omega
_{t}\left( G-1\right) \geq 1$.

$ii)$ By Proposition 5.3 (with $\psi $ instead of $q)$, $D_{\psi }\left(
G\right) =G.D_{\psi }F$, i.e., by the definition of $D_{\psi }$: 
\begin{equation*}
\left( G(\psi t)-G(t)\right) /\left( \psi -1\right) t=G(t).D_{\psi }F(t)%
\text{,}
\end{equation*}
therefore: 
\begin{equation}
G(t)=\left( 1-\left( 1-\psi \right) tD_{\psi }F(t\right) )^{-1}G\left( \psi
t\right) \text{.}  \tag{8.2}
\end{equation}

Iterating $\left( 8.2\right) $ gives, for all integers $k\geq 1$, $%
G(t)=\prod\nolimits_{j=0}^{k-1}\left( 1-\left( 1-\psi \right) \psi
^{j}t\left( D_{\psi }F\right) \left( \psi ^{j}t\right) \right) ^{-1}G\left(
\psi ^{k}t\right) $, where all the fps of this equality belong to $B\left[ %
\left[ q\right] \right] \left[ \left[ t\right] \right] =B\left[ \left[ q,t%
\right] \right] $. From $i)$ we have $G\left( \psi ^{k}t\right)
-1=\sum\nolimits_{n\geq 1}\gamma _{n}\left( \psi ^{k}t\right) ^{n}/\left[ n%
\right] _{\psi }!$, and by hypothesis, $\omega (\psi )$ $\geq 1$. Therefore, 
$\omega \left( G\left( \psi ^{k}t\right) -1\right) \geq \omega _{t}\left(
G-1\right) .\left( k+1\right) \geq k+1$, so $G\left( \psi ^{k}t\right) $
converges formally to $1$ in $B\left[ \left[ q,t\right] \right] $, which
implies the formal convergence of the infinite product.
\end{proof}

\bigskip Note that for $F\left( t\right) =t$ and $\psi \left( q\right) =q$,
Eq. $\left( 8,1\right) $ gives the well-known formula $\mathbf{e}_{q}\left(
t\right) =\prod\nolimits_{k\geq 0}(1-\left( 1-q\right) q^{k}t)^{-1}$.

\begin{corollary}
In $\Lambda _{\mathbb{Q}}\left[ \left[ q,t\right] \right] $, we have the
factorisation 
\begin{equation}
E(t)=\sum\limits_{n=0}^{\infty }e_{n}t^{n}=\prod\limits_{k=0}^{\infty
}\left( 1-\left( 1-q\right) \left( q^{k}t\right) p_{q}\left( -q^{k}t\right)
\right) ^{-1}=\prod\limits_{k=0}^{\infty }\left( 1+\left( 1-q\right)
\sum\nolimits_{n\geq 1}\left[ p_{n}\right] _{q}\left( -q^{k}t\right)
^{n}\right) ^{-1}  \tag{8.3}
\end{equation}
where the series and product are formally convergent.
\end{corollary}

\begin{proof}
\bigskip Take $B=\mathbb{Q}\left[ \left[ X\right] \right] $ (with $X=\left(
x_{i}\right) $), $F\left( t\right) =-P_{q}(-t)$, and $\psi \left( q\right)
=q $. Definition 8.1 gives $D_{q}F=D_{q}\left( -P_{q}\left( -t\right)
\right) =D_{q}\left( P_{q}\right) \left( -t\right) =p_{q}\left( -t\right)
=\sum\nolimits_{n\geq 1}\left[ p_{n}\right] \left( -t\right) ^{n-1}$. Hence,
using $\left( 8.1\right) $ and Theoreme 5.5, we obtain Eq. $\left(
8.3\right) $. Furthermore, $\Lambda _{\mathbb{Q}}$ is a subalgebra of $B$;
hence, the infinite product is in $\Lambda _{Q}\left[ \left[ q,t\right] %
\right] $.
\end{proof}

\bigskip

\subsection{\protect\bigskip Case of $H(t)$}

Since $H(t)=\left( E(-t)\right) ^{-1}$, its infinite product follows
directly from $\left( 8.3\right) $. However, we will find this product by
another method.

\begin{proposition}
\bigskip Let $F\in B\left[ \left[ q,t\right] \right] $ with $\omega
_{t}\left( F\right) \geq 1$,$\ $and let $\psi \left( q\right) =q^{m}$ with $%
m\in \mathbb{Z}^{\ast }$.

$i)$ The fps $G=\mathbf{E}_{\psi }\left[ F\right] _{\psi }^{\ast }$ belongs
to $B\left[ \left[ q,t\right] \right] $.

$ii)$ We have 
\begin{equation}
\mathbf{E}_{\psi }\left[ F\right] _{\psi }^{\ast
}=\prod\limits_{k=0}^{\infty }\left( 1+\left( 1-\psi \right) \psi
^{k}t\left( D_{\psi }F\right) \left( \psi ^{k}t\right) \right) \text{,} 
\tag{8.4}
\end{equation}
where the convergence of the product is that of the formal convergence in $B%
\left[ \left[ q,t\right] \right] $.
\end{proposition}

\begin{proof}
$i)$ Let $F(t)=\sum\nolimits_{n\geq 1}f_{n}t^{n}/\left[ n\right] _{\psi }!$.
As before, $B\left[ \left[ q\right] \right] \subset \mathbb{Q}\left( \left(
q\right) \right) \left[ \left[ Y\right] \right] $, and according to
Proposition 6.6 and Section 7 with $A=\mathbb{Q}\left( \left( q\right)
\right) \left[ \left[ Y\right] \right] $, the fps $\mathbf{E}_{\psi }\left[ F%
\right] _{\psi }^{\ast }=\sum\nolimits_{n=0}^{\infty }\gamma _{n}^{\ast
}\left( t^{n}/\left[ n\right] _{\psi }!\right) $ is well-defined in $A\left[ %
\left[ t\right] \right] $, and by $\left( 6.6\right) $, $\gamma _{0}^{\ast
}=1$. For $n\geq 0$,$\ \ \gamma _{n+1}^{\ast }=\sum\nolimits_{k=0}^{n}%
\QATOPD[ ] {n}{k}_{\psi }\psi ^{n-k}\gamma _{n-k}^{\ast }f_{k+1}$. By the
assumption on $F$ and $\psi $, we have $f_{n}/\left( \left[ n\right] _{\psi
}!\right) $ $\in B\left[ \left[ q\right] \right] $, and $\left[ n\right]
_{\psi }!\in \mathbb{Q}\left[ \left[ q\right] \right] $, therefore $f_{n}\in
B\left[ \left[ q\right] \right] $. Furthermore, $\QATOPD[ ] {n}{j}_{\psi }$%
and $\psi ^{n-k}$ are polynomials in $\psi $, so we deduce by induction that
the coefficients $\gamma _{n}^{\ast }$ and also $\gamma _{n}^{\ast }/\left[ n%
\right] _{\psi }!$ belong to $B\left[ \left[ q\right] \right] $. Thus, $%
G^{\ast }\left( t\right) =\mathbf{E}_{\psi }\left[ F\right] _{\psi }^{\ast
}\in B\left[ \left[ q,t\right] \right] $.

$ii)$ By Proposition 6.4 (with $\psi $ instead of $q)$, \ $D\psi \left(
G^{\ast }\right) =G^{\ast }\left( \psi t\right) .D_{\psi }F$. So, by the
definition of $D_{\psi }$, $G^{\ast }(\psi t)-G^{\ast }(t)/\left( \psi
-1\right) t=G^{\ast }(\psi t).D_{\psi }F(t)$, and we get 
\begin{equation*}
G^{\ast }(t)=1+\left( 1-\psi \right) tD_{\psi }F\left( t\right) )G^{\ast
}\left( \psi t\right) \text{.}
\end{equation*}

The rest of the proof is analogous to the previous case.
\end{proof}

Note that for $F\left( t\right) =t$ and $\psi \left( q\right) =q$, Eq. $%
\left( 8.4\right) $ gives the well-known formula $\mathbf{E}_{q}\left(
t\right) =\prod\nolimits_{k\geq 0}\left( 1+\left( 1-q\right) q^{k}t\right) .$

\bigskip As before we deduce the following:

\begin{corollary}
\bigskip In $\Lambda _{\mathbb{Q}}\left[ \left[ q,t\right] \right] $ we have
the factorization 
\begin{equation}
H(t)=\sum\limits_{n=0}^{\infty }h_{n}t^{n}=\prod\limits_{k=0}^{\infty
}\left( 1+\left( 1-q\right) \left( q^{k}t\right) p_{q}\left( q^{k}t\right)
\right) =\prod\limits_{k=0}^{\infty }\left( 1+\left( 1-q\right)
\sum\nolimits_{n\geq 1}\left[ p_{n}\right] \left( q^{k}t\right) ^{n}\right) 
\text{,}  \tag{8.5}
\end{equation}
where the series and products are formally convergent.\bigskip
\end{corollary}

\bigskip This is, of course, what we obtain directly with $\left( 8.3\right) 
$ and $H\left( t\right) =\left( E\left( -t\right) \right) ^{-1}$.

\begin{proposition}
\bigskip Let $F\in B\left[ \left[ q,t\right] \right] $ with $\omega
_{t}\left( F\right) \geq 1$, then 
\begin{equation}
\mathbf{e}_{q}\left[ -F(t)\right] _{q}=\left( \mathbf{E}_{q}\left[ F(t)%
\right] _{q}^{\ast }\right) ^{-1}\text{,}  \tag{8.6}
\end{equation}
or equivalently 
\begin{equation}
\mathbf{E}_{q}\left[ -F(t)\right] _{q}^{\ast }=\left( \mathbf{e}_{q}\left[
F(t)\right] _{q}\right) ^{-1}\text{.}  \tag{8.7}
\end{equation}
\end{proposition}

We observe that $\left( 8.6\right) $ and $\left( 8.7\right) $ are two
analogs of the well-known following equation for the classical formal
composition of fps with $F(0)=0$: $\exp \left( -F(t)\right) =\left( \exp
F(t)\right) ^{-1}$.

\begin{proof}
Starting from $\left( 8.1\right) $ with $\psi \left( q\right) =q$%
\begin{equation*}
\mathbf{e}_{q}\left[ -F(t)\right] _{q}=\prod\limits_{k=0}^{\infty }\left(
1-\left( 1-q\right) \right) \left( q^{k}t\right) D_{q}(-F)\left(
q^{k}t\right) )^{-1}\text{.}
\end{equation*}

It is easy to check that $D_{q}\left( -F\right) \left( t\right) =-D_{q}F(t)$%
. From this, we have 
\begin{equation*}
\mathbf{e}_{q}\left[ -F\right] _{q}=\prod\limits_{k=0}^{\infty }\left(
1+\left( 1-q\right) \right) \left( q^{k}t\right) D_{q}F\left( q^{k}t\right)
)^{-1}\text{.}
\end{equation*}
In comparing this last equation with $\left( 8.4\right) $ and $\psi \left(
q\right) =q$, we obtain $\left( 8.6\right) $. Eq. $\left( 8.7\right) $ is
obtained by replacing $F$ by $-F$.
\end{proof}

\section{\protect\bigskip Applications by specialization}

\bigskip

\subsection{\protect\bigskip General case}

The specialization of symmetric functions is a crucial concept for
understanding this section. Recall that this involves a unitary homomorphism
between $\Lambda _{K}$ and a commutative algebra over $K$. This homomorphism
can be defined either by its value on the variables $x_{i}$ or by the values
assigned to $e_{n}$ or $h_{n}$, as these functions are algebraically
independant and generate $\Lambda _{K}$. A detailed description of
specialization can be found in [15, Chap.7].

Consider a formal series: 
\begin{equation}
G(t)=1+\sum\nolimits_{n\geq 1}a_{n}t^{n}\text{.}  \tag{9.1}
\end{equation}
We can apply the results from the previous sections by considering $G(t)$ as
a specialization of the generating series of symmetric functions. This can
be done in two ways: either as the specialization of $E(t)$ or as the
specialization of $H(t)$. In the first case the homomorphisme sends $e_{n}$
to $a_{n}$, and in the second case, it sends $h_{n}$ to $a_{n}$. We will
apply the results of Sections 3 to 7 with $a_{n}\in A=\mathbb{C}\left(
\left( q\right) \right) \left[ \left[ Y\right] \right] $, where $Y$ is a
finite or infinite sequence of variables (distinct, of course, from $q$ and $%
t$). For the application of Section 8, we will limit ourselves to cases
where $a_{n}\in A=\mathbb{C}\left[ \left[ q,Y\right] \right] $. The
convergence of infinite series and products is always that of formal
convergence.

\bigskip

\textbf{Specialization }$E(t)\equiv G(t)$

By identification with $\left( 3.1\right) $, we have 
\begin{equation}
p_{q}\left( -t\right) =\dfrac{D_{q}G(t)}{G(t)}\text{,}  \tag{9.2}
\end{equation}
\ where 
\begin{equation}
\;p_{q}(t)=D_{q}P_{q}(t)=\sum\limits_{n\geq 0}\left[ p_{n+1}\right] t^{n}%
\text{.}  \tag{9.3}
\end{equation}
The $\left[ p_{n}\right] $ are related to $a_{n}$ by Eq. $\left( 3.2\right) $%
, $\left( 3.4\right) $ and $\left( 3.5\right) $ with $e_{n}$ replaced by $%
a_{n}$. Sections 5 and 8 \ give 
\begin{equation}
G\left( t\right) =\mathbf{e}_{q}\left[ -P_{q}(-t)\right] _{q}\text{,} 
\tag{9.4}
\end{equation}
and 
\begin{equation}
G(t)=\prod\limits_{k=0}^{\infty }\left( 1-\left( 1-q\right)
q^{k}tp_{q}(-q^{k}t)\right) ^{-1}=\prod\limits_{k=0}^{\infty }\left(
1+\left( 1-q\right) \sum\nolimits_{n\geq 1}\left[ p_{n}\right] \left(
-q^{k}t\right) ^{n}\right) ^{-1}\text{.}  \tag{9.5}
\end{equation}

\bigskip

\textbf{Specialization }$H(t)\equiv G(t)$

By identification with $\left( 3.6\right) $ and introducing the exponent $h$
to distinguish these expressions from those defined in the previous case, we
have: 
\begin{equation}
p_{q}^{h}(t)=\dfrac{D_{q}G(t)}{G\left( qt\right) }\text{,}  \tag{9.6}
\end{equation}
where $p_{q}^{h}(t)=D_{q}P_{q}^{h}(t)=\sum\limits_{n\geq 0}\left[ p_{n+1}^{h}%
\right] t^{n}$. The $\left[ p_{n}^{h}\right] $ are given this time by Eq. $%
\left( 3.7\right) $, $\left( 3.8\right) $ and $\left( 3.9\right) $, with $%
h_{n}$ replaced by $a_{n}$. Sections 6 and 8 give\ $G\left( t\right) =%
\mathbf{E}_{q}\left[ P_{q}^{h}(t)\right] _{q}^{\ast }$ and

\begin{equation}
G(t)=\prod\limits_{k=0}^{\infty }\left( 1+\left( 1-q\right)
q^{k}tp_{q}^{h}(q^{k}t)\right) =\prod\limits_{k=0}^{\infty }\left( 1+\left(
1-q\right) \sum\nolimits_{n\geq 1}\left[ p_{n}^{h}\right] (q^{k}t)^{n}%
\right) \text{.}  \tag{9.7}
\end{equation}
It follows from the above, that any non-zero formal power series with a
constant term equal to 1 (which can always be achieved through
normalization), can be expressed, using the $\left[ p_{n}\right] $, in the
form of a $q$-exponential formula or an infinite $q$-product. However, for
these expressions not to yield trivial results, certain conditions on $G(t)$
are necessary. For example, $G(t)$ should not be equal to a monomial.
Furthermore, for the infinite $q$-product to be interesting, it will
necessary to be able to calculate the formal sums of the fps appearing in
each of the factors of the product. We will examine some examples.

\subsection{$q$-binomial Theorem}

In this elementary example, we show that $q$-binomial theorem is a
particular case of the infinite $q$-product formula $\left( 9.5\right) $. We
refer to [7] for the usual proof of this theorem and adopt the following
standard notations:

i) If $a=0$, \ $\left( a;q\right) _{n}=1$ \ \ for $n\in \mathbb{N}$ or $%
n=\infty $

\bigskip ii) If $a\neq 0$%
\begin{equation*}
\left( a;s\right) _{n}\left\{ 
\begin{array}{c}
=1\;\;\;\;\;\;\text{if }n=0\;\;\;\ \ \ \ \ \ \ \ \ \ \ \ \ \ \ \ \ \ \ \ \ \
\ \ \ \ \ \ \;\;\;\;\;\;\;\;\; \\ 
=\left( 1-a\right) \left( 1-qa\right) ...\left( 1-q^{n-1}a\right) \text{ \
if }n=1,2,...
\end{array}
\right.
\end{equation*}
and $\left( a;q\right) _{\infty }=\prod\nolimits_{k\geq 0}\left(
1-q^{k}a\right) $.

Let

\begin{equation*}
G(t)=\sum\limits_{n=0}^{\infty }\dfrac{\left( a;q\right) _{n}}{\left(
q,q\right) _{n}}t^{n}\text{,}
\end{equation*}
and consider $G(t)$ as the specialization of $E(t)$. Using $\left(
3.4\right) $, we have: 
\begin{equation*}
\left[ p_{1}\right] =\left| \left[ 1\right] a_{1}\right| =\dfrac{1-a}{1-q}%
\text{, }\left[ p_{2}\right] =\left| 
\begin{tabular}{ll}
$\left[ 1\right] a_{1}$ & $\left[ 1\right] $ \\ 
$\left[ 2\right] a_{2}$ & $a_{1}$%
\end{tabular}
\right| =\left( -a\right) \dfrac{1-a}{1-q}\text{,}
\end{equation*}
and by an immediate recurrence with $\left( 3.3\right) $: 
\begin{equation*}
\left[ p_{n}\right] _{q}=\left( -a\right) ^{n-1}\dfrac{1-a}{1-q}\text{,}
\end{equation*}
from which, by Definition 8.1 $ii)$ 
\begin{equation*}
p_{q}(t)=\dfrac{1-a}{1-q}\sum\limits_{n\geq 0}\left( -at\right) ^{n}=\dfrac{%
1-a}{\left( 1-q\right) \left( 1+at\right) }\text{,}
\end{equation*}
and by $\left( 5.8\right) $%
\begin{equation*}
P_{q}\left( t\right) =\dfrac{a-1}{a}\sum\limits_{n\geq 1}\dfrac{\left(
-at\right) ^{n}}{1-q^{n}}\text{.}
\end{equation*}
Therefore, by Theorem 5.5, a new formula is 
\begin{equation}
G(t)=\mathbf{e}_{q}\left[ \dfrac{1-a}{a}\sum\limits_{n\geq 1}\dfrac{\left(
at\right) ^{n}}{1-q^{n}}\right] _{q}\text{,}  \tag{9.8}
\end{equation}
and $\left( 9.5\right) $ gives: 
\begin{equation*}
G(t)=\prod\limits_{k=0}^{\infty }\left( 1+\left( 1-q\right)
\sum\nolimits_{n\geq 1}\left( -a\right) ^{n-1}\dfrac{1-a}{1-q}\left(
-q^{k}t\right) ^{n}\right) ^{-1}=\prod\limits_{k=0}^{\infty }\left( 1-\left(
1-a\right) \dfrac{q^{k}t}{1-aq^{k}t}\right) ^{-1}\text{.}
\end{equation*}
After simplifying each of the factors in the last product above, we obtain
the $q$-binomial theorem: 
\begin{equation*}
\sum\limits_{n=0}^{\infty }\dfrac{\left( a;q\right) _{n}}{\left( q,q\right)
_{n}}t^{n}=\dfrac{\left( at;q\right) _{\infty }}{\left( t;q\right) _{\infty }%
}\text{.}
\end{equation*}

\subsection{Inversions in trees}

Consider the $q$-deformation of the exponential series $E_{xp}(t)=\sum%
\nolimits_{n\geq 0}q^{\binom{n}{2}}t^{n}/\left[ n\right] !$, encountered in
[4] (and also in [5]). We specialize by setting $E(t)\equiv E_{xp}(t)$. It
was shown in [4, Eq. $\left( 6.4\right) $ for $r=1$] that $p_{n}=\left(
1-q\right) ^{n-1}J_{n}(q)/\left( n-1\right) !$, and in [4, Lemma 6.1 pour $%
r=1$], that \ $p_{n}(q)=\left( 1-q\right) \left[ p_{n-1}\right] _{q}$.
Therfore, we have:

\begin{equation*}
\left[ p_{n}\right] =\dfrac{\left( 1-q\right) ^{n-1}}{n!}J_{n+1}(q)\text{.}
\end{equation*}
Recall that $J_{n}\left( q\right) =\sum\nolimits_{T}q^{inv(T)}$, where the
sum is over the set of trees on vertices $\mathbf{n}$, and where $inv(T)$ is
the number of inversions of $T$. Recall also that an inversion of a tree $T$
is a pair of vertices $\left( i,j\right) $, such that, $i>j>1$, and such
that the unique path from the root $1$ to $j$ passes through $i$ (see [14]
or [15, Chap.5] for more details).

After simplification, Formulas $\left( 3.2\right) ,\left( 3.4\right) $ and\ $%
\left( 3.5\right) $ give successively:

\begin{equation}
\left[ n\right] q^{\binom{n}{2}}=\sum\limits_{k=1}^{n}\binom{n}{k}q^{\binom{%
n-k}{2}}\left( q-1\right) ^{k-1}J_{k+1}(q)\text{,}  \tag{9.9}
\end{equation}
$\qquad \qquad \qquad \qquad \qquad \qquad $%
\begin{equation}
\left( 1-q\right) ^{n-1}J_{n+1}(q)=n!\left| 
\begin{array}{ccccccc}
\left[ 1\right] & 1 & 0 & 0 & . & . & 0 \\ 
\left[ 2\right] q/2! & 1 & 1 & 0 & . & . & 0 \\ 
. & . & . & . & . & . & . \\ 
. & . & . & . & . & . & . \\ 
. & . & . &  & . & . & 0 \\ 
. & . & . &  &  & . & 1 \\ 
\left[ n\right] q^{\binom{n}{2}}/n! & q^{\binom{n-1}{2}}/(n-1)! & . & . & .
& q/2! & 1
\end{array}
\right|  \tag{9.10}
\end{equation}
$\qquad \qquad \qquad \qquad \qquad \qquad \qquad \qquad $%
\begin{equation}
\dfrac{\left[ n\right] !}{n!}q^{\binom{n}{2}}=\left| 
\begin{array}{ccccccc}
J_{2} & \left[ 1\right] & 0 & 0 & . & . & 0 \\ 
\left( 1-q\right) \dfrac{J_{3}}{2!} & J_{2} & \left[ 2\right] & 0 & . & . & 0
\\ 
. & . & . & . & . & . & . \\ 
. & . & . & . & . & . & . \\ 
. & . &  & . & . & . & 0 \\ 
. & . &  &  &  & . & \left[ n-1\right] \\ 
\left( 1-q\right) ^{n-1}\dfrac{J_{n+1}}{n!} & \left( 1-q\right) ^{n-1}\dfrac{%
J_{n}}{n-1!} & . & . & . & \left( 1-q\right) \dfrac{J_{3}}{2!} & J_{2}
\end{array}
\right|  \tag{9.11}
\end{equation}
From (9.8) and (8.3), we obtain 
\begin{equation}
E_{xp}(t)=\sum\nolimits_{n\geq 0}q^{\binom{n}{2}}\dfrac{t^{n}}{n!}=\mathbf{e}%
_{q}\left[ \sum\limits_{n\geq 1}\dfrac{\left( q-1\right) ^{n-1}}{\left[ n%
\right] n!}J_{n+1}(q)t^{n}\right] _{q}\text{,}  \tag{9.12}
\end{equation}

\begin{equation}
E_{xp}(t)=\prod\limits_{k\geq 0}\left( 1+\sum\nolimits_{n\geq 1}\dfrac{%
\left( 1-q\right) ^{n}}{n!}\left( -q^{k}t\right) ^{n}J_{n+1}(q).\right) ^{-1}%
\text{.}  \tag{9.13}
\end{equation}
In this last formula, it does not appear possible to simplify the series in
each factor, which limits its practical interest. Note that when $q=1$, Eq. $%
\left( 9.9\right) $ gives a trivial result and Eq. $\left( 9.10\right) $ to $%
\left( 9.13\right) $ have no counterparts. Moreover, except for $\left(
9.9\right) $ which we discuss later, these formulas seem to be new.

The reciprocal polynomials of $J_{n}$ have also been defined in [4] as: 
\begin{equation}
\overline{J_{n}}(q)=q^{\binom{n-1}{2}}J_{n}\left( 1/q\right) \Leftrightarrow
J_{n}(q)=q^{\binom{n-1}{2}}\overline{J_{n}(q)\text{.}}  \tag{9.14}
\end{equation}
\qquad These reciprocal polynomials are equal to the sum enumerator of the
parking function (see [16] for more details). With $\left( 9.14\right) $, we
can obtain the equations corresponding to $\left( 9.9\right) $, $\left(
9.10\right) $ and $\left( 9.11\right) $ for reciprocal polynomials. Thus, $%
\left( 9.9\right) $ yields:

\begin{equation}
\left[ n\right] =\sum\limits_{k=1}^{n}\binom{n}{k}q^{\left( k+1\right)
\left( n-k\right) }\left( q-1\right) ^{k-1}\overline{J_{k+1}(q)}  \tag{9.15}
\end{equation}
This equation is equivalent to a particular case of the linear equation $%
\left( 1.27\right) $ from [12, p. 37]. This equation, which we have already
encountered in [4], was derived in [12] by applying Goncarov's polynomial
theory to parking functions. When $u_{i}=i$ and taking into account that $%
S_{k}\left( q\right) =\overline{J_{k+1}(q)}$, Eq. $\left( 1.27\right) $
becomes for the specific case of interest: 
\begin{equation}
1=\sum\limits_{k=0}^{n}\binom{n}{k}q^{\left( k+1\right) \left( n-k\right)
}\left( 1-q\right) ^{k}\overline{J_{k+1}(q)}\text{,}  \tag{9.16}
\end{equation}
It is easy to see the equivalence of $\left( 9.16\right) $ with $\left(
9.15\right) $. Additionally, we could obtain with $\left( 9.14\right) $ the
equation corresponding to $\left( 9.12\right) $ for reciprocal polynomials.
However, we can find this corresponding equation in a more general manner,
as Eq. $\left( 9.12\right) $\ is equivalent to Recurrence $\left( 5.7\right) 
$, whith $\gamma _{n}=q^{\binom{n}{2}}\left[ n\right] !/n!$ and $%
f_{n}=\left( q-1\right) ^{n-1}(\left[ n-1\right] !/n!)J_{n+1}$. The
calculation we performed in the proof b) of Theorem 6.7 is actually valid
for any pair of sequences $\left( \gamma _{n},f_{n}\right) $ satisfying
(5.7), by setting: 
\begin{equation*}
\overline{\gamma _{n}}=q^{\binom{n}{2}}\gamma _{n}\left( 1/q\right) =q^{-%
\binom{n}{2}}\dfrac{\left[ n\right] !}{n!}\text{ \ \ \ et \ \ \ }\overline{%
f_{n}}=q^{\binom{n}{2}}f_{n}\left( 1/q\right) =q^{-\binom{n}{2}}\left(
1-q\right) ^{n-1}\dfrac{\left[ n-1\right] !}{n!}\overline{J_{n+1}}\text{.\ \ 
}
\end{equation*}
Replacing these into $\left( 6.9\right) $ yields the corresponding
sought-after equation:

\begin{equation*}
\sum\nolimits_{n\geq 0}q^{-\binom{n}{2}}\dfrac{t^{n}}{n!}=\mathbf{E}_{q}%
\left[ \sum\limits_{n\geq 1}q^{-\binom{n}{2}}\dfrac{\left( 1-q\right) ^{n-1}%
}{\left[ n\right] n!}\overline{J_{n+1}}(q)t^{n}\right] _{q}^{\ast }\text{,}
\end{equation*}
where the summation is to be taken in $\mathbb{Q}\left( \left( q\right)
\right) \left[ \left[ t\right] \right] $.

\subsubsection{\protect\bigskip $q$-orthogonal polynomials}

In this example, we apply the results of the previous sections in the
reverse manner compared to the first example: it is now the knowledge of the
infinite product in $\left( 9.4\right) $ that will allow us to deduce other
formulas for these polynomials. We refer to [9] and [11] for more details on
the definition and properties of these $q$-orthogonal polynomials, and use,
as in those references, the variables $x$ and $q$ for these polynomials. It
is important not to confuse this $x$ -which plays the role of $Y$ in the
general case presented in Subsection 9.1- with the variables $X=\left(
x_{i}\right) $ of symmetric functions. We provide a detailed treatment of
the case of discrete $q$-Hermite polynomials. We denote the discrete $q$%
-Hermite I polynomials by an uppercase letter, $H_{n}\left( x;q\right) $, to
avoid any confusion with complete symmetric functions $h_{n}$. Let $G$ be
the generating function of discrete $q$-Hermite I polynomials given in [11,
Eq. (14.28.11)] by:

\begin{equation}
G(t,x)=\sum\limits_{n=0}^{\infty }\dfrac{H_{n}\left( x;q\right) }{\left(
q,q\right) _{n}}t^{n}=\dfrac{\left( t^{2};q^{2}\right) _{\infty }}{\left(
xt;q\right) _{\infty }}\text{.}  \tag{9.17}
\end{equation}
Equation $\left( 9.9\right) $ expands as follows:$\qquad \qquad \qquad $%
\begin{equation}
G(t,x)=\prod\limits_{k=0}^{\infty }\left( \dfrac{1-xtq^{k}}{1-t^{2}q^{2k}}%
\right) ^{-1}\text{.}  \tag{9.18}
\end{equation}
Let's specialize $E(t)=G(t,x)$. By comparing the factors of $\left(
9.17\right) $ and $\left( 8.3\right) $ corresponding to $k=0$, we are led to
assume 
\begin{equation*}
\dfrac{1-xt}{1-t^{2}}=1-\left( 1-q\right) tp_{q}\left( -t\right) \text{,}
\end{equation*}
or equivalently 
\begin{equation}
p_{q}(q)=\dfrac{x+t}{\left( 1-q\right) \left( 1-t^{2}\right) }\text{.} 
\tag{9.19}
\end{equation}
Let's now prove $\left( 9.19\right) $. Set $f\left( t\right) =\left(
x-t\right) /(\left( 1-q)(1-t^{2}\right) )$. Expanding $f(t)$, we obtain: 
\begin{equation}
f(t)=\left( 1-q\right) ^{-1}\left\{
x-t+xt^{2}-t^{3}+xt^{4}-t^{5}+...\right\} \text{.}  \tag{9.20}
\end{equation}
Let $F(t)$ be the fps whose $q$-derivative is $f(t)$, and which has a zero
constant term, that is: 
\begin{equation}
F\left( t\right) =\left( 1-q\right) ^{-1}\left\{ xt-\dfrac{t^{2}}{\left[ 2%
\right] }+x\dfrac{t^{3}}{\left[ 3\right] }-\dfrac{t^{4}}{\left[ 4\right] }%
+...\right\} \text{.}  \tag{9.21}
\end{equation}
According to Proposition 8.2 with $\psi \left( q\right) =q$, we have 
\begin{equation*}
\mathbf{e}_{q}\left[ F\left( t\right) \right] _{q}=\prod\limits_{k\geq
0}\left( 1-\left( 1-q\right) q^{k}tf\left( q^{k}t\right) \right)
=\prod\limits_{k=0}^{\infty }\left( \dfrac{1-xtq^{k}}{1-t^{2}q^{2k}}\right)
^{-1}\text{.}
\end{equation*}
Therefore, with $\left( 9.4\right) $ and (9.18): 
\begin{equation*}
\mathbf{e}_{q}\left[ F\left( t\right) \right] =G(t,x)=\mathbf{e}_{q}\left[
-P_{q}\left( -t\right) \right] \text{.}
\end{equation*}
Now, the mapping that associates the $q$-composition by $\mathbf{e}_{q\text{ 
}}$to a formal series is injective, as shown by the inversion of Eq. $\left(
5.7\right) $, given by: 
\begin{equation*}
f_{0}=0\text{, and for }n\geq 0\text{,}\;\;f_{n+1}=\gamma
_{n+1}-\sum\limits_{k=1}^{n-1}\QATOPD[ ] {n}{k}\gamma _{n-k}f_{k+1}\text{.}
\end{equation*}
Thus, we have indeed verified that $F\left( t\right) =-P_{q}\left( -t\right) 
$ and by $q$-derivation $\left( 9.19\right) $. Therefore, using Eq. $\left(
5.9\right) $, we obtain the generating function in the form of a $q$%
-composition:

\begin{equation}
G(t,x)=\sum\limits_{n=0}^{\infty }\dfrac{H_{n}\left( x;q\right) }{\left(
q;q\right) _{n}}t^{n}=\mathbf{e}_{q}\left[ \left( 1-q\right) ^{-1}\left\{ xt-%
\dfrac{t^{2}}{\left[ 2\right] }+x\dfrac{t^{3}}{\left[ 3\right] }-\dfrac{t^{4}%
}{\left[ 4\right] }+...\right\} \right] _{q}\text{.}  \tag{9.22}
\end{equation}
This generating function is, modulo variable changes, a $q$-analog of the
following generating function for standard Hermite polynomials $H_{n}\left(
x\right) $, [11, Eq. (9.15.10)]:

\begin{equation*}
\sum\limits_{n=0}^{\infty }\dfrac{H_{n}\left( x\right) }{n!}t^{n}=\exp
\left( 2x-t^{2}\right) \text{.}
\end{equation*}
Indeed, we know that\ [11, Eq. 14.28.12]:

\begin{equation*}
\lim\limits_{q\rightarrow 1}\dfrac{H_{n}\left( x\sqrt{1-q^{2}};q\right) }{%
\left( 1-q^{2}\right) ^{n/2}}=\dfrac{H_{n}\left( x\right) }{2^{n}}\text{,}
\end{equation*}
thus, 
\begin{equation*}
\dfrac{H_{n}\left( x\right) }{2^{n}n!}=\lim\limits_{q\rightarrow 1}\dfrac{%
H_{n}\left( x\sqrt{1-q^{2}};q\right) }{\left( 1-q^{2}\right) ^{n/2}}\dfrac{%
\left( 1-q\right) ^{n}}{\left( q;q\right) _{n}}=\lim\limits_{q\rightarrow 1}%
\dfrac{H_{n}\left( x\sqrt{1-q^{2}};q\right) }{\left( q;q\right) _{n}}\left( 
\dfrac{1-q}{1+q}\right) ^{n/2}\text{,}
\end{equation*}
and consequently, with $\mu =\left( (1-q)/\left( 1+q\right) \right) ^{1/2}$: 
\begin{equation*}
\lim\limits_{q\rightarrow 1}\mathbf{e}_{q}\left[ \left( 1-q\right)
^{-1}\left\{ x\sqrt{1-q^{2}}\mu t-\dfrac{\left( \mu t\right) ^{2}}{\left[ 2%
\right] }+x\sqrt{1-q^{2}}\dfrac{\left( \mu t\right) ^{3}}{\left[ 3\right] }-%
\dfrac{\left( \mu t\right) ^{4}}{\left[ 4\right] }+...\right\} \right]
_{q}=\exp \left( xt-\dfrac{t^{2}}{4}\right) \text{,}
\end{equation*}
where the limit must be understood in the sense of convergence with respect
to $q$, and for each $n$, from the coefficient of $t^{n}$ in the series on
the left-hand side to the coefficient of $t^{n}$ in the series on the
right-hand side. Furthermore, according to Definition 8.1 $ii)$, Eq. $\left(
9.20\right) $, and $p_{q}\left( t\right) =f(-t)$, we have:

\begin{equation}
\left\{ 
\begin{array}{c}
\text{For odd }n\text{, }\left[ p_{n}\right] =x\left( 1-q\right) ^{-1} \\ 
\text{For even }n\text{, }\left[ p_{n}\right] =\left( 1-q\right) ^{-1}
\end{array}
\right.  \tag{9.23}
\end{equation}
Therefore, using $\left( 3.5\right) $, we obtain the following determinants
for the discrete $q$-Hermite I polynomials:

\begin{equation}
\text{For }n\text{ odd, }H_{n}\left( x;q\right) =\left| 
\begin{array}{ccccccc}
x & 1-q & 0 & 0 & . & . & 0 \\ 
1 & x & 1-q^{2} & 0 & . & . & 0 \\ 
x & 1 & x & 1-q^{3} & . & . & . \\ 
. & . & . & . & . & . & . \\ 
. & . & . &  & . & . & 0 \\ 
. & . & . &  &  & . & 1-q^{n-1} \\ 
x & 1 & x & . & . & 1 & x
\end{array}
\right|  \tag{9.24}
\end{equation}
$\qquad \qquad \qquad \qquad \qquad $%
\begin{equation}
\text{For }n\text{ even, }H_{n}\left( x;q\right) =\left| 
\begin{array}{ccccccc}
x & 1-q & 0 & 0 & . & . & 0 \\ 
1 & x & 1-q^{2} & 0 & . & . & 0 \\ 
x & 1 & x & 1-q^{3} & . & . & . \\ 
. & . & . & . & . & . & . \\ 
. & . & . &  & . & . & 0 \\ 
. & . & . &  &  & . & 1-q^{n-1} \\ 
1 & x & 1 & . & . & 1 & x
\end{array}
\right|  \tag{9.25}
\end{equation}
Now, specializing $H(t)=G(t,x)$, we have 
\begin{equation*}
p_{q}^{h}(t)=\dfrac{D_{q}G(t,x)}{G(qt,x)}=p_{q}(-t)\dfrac{G(t,x)}{G(qt,x)}=%
\dfrac{x-t}{\left( 1-q\right) \left( 1-t^{2}\right) }\dfrac{\left(
t^{2};q^{2}\right) _{\infty }}{\left( xt;q\right) _{\infty }}\dfrac{\left(
xqt;q\right) _{\infty }}{\left( q^{2}t^{2};q\right) _{\infty }}=\dfrac{x-t}{%
\left( 1-q\right) \left( 1-xt\right) }\text{.}
\end{equation*}
Expanding, we obtain 
\begin{equation*}
p_{q}^{h}(t)=\left( 1-q\right) ^{-1}\left\{ x+\left( x^{2}-1\right) t+\left(
x^{3}-x\right) t^{2}+...\right\} \text{,}
\end{equation*}
which implies 
\begin{equation*}
\left\{ 
\begin{array}{c}
\text{ }\left[ p_{1}^{h}\right] =\left( 1-q\right) ^{-1}x \\ 
\text{Pour }n\geq 2\text{ \ }\left[ p_{n}\right] =\left( 1-q\right)
^{-1}x^{n-2}\left( x^{2}-1\right)
\end{array}
\right. \text{,}
\end{equation*}
and with $\left( 3.9\right) $, we find for $n\geq 2$ and after
simplification, the formula $\left( 1.6\right) $ from the introduction.To
the best of the author's knowledge, the formulas $\left( 9.22\right) $, $%
\left( 9.24\right) $, $\left( 9.25\right) $ and $\left( 1.6\right) $ are new$%
^{b}$ 
\footnotetext{$^{b}$\textit{However, it is possible that I may have missed
something in the literature on this topic. I invite the reader who is aware
of relevant references to share them with me. }}.

\begin{remark}
It is known that there exist representations of orthogonal polynomials in
the form of moment determinants (see, for example, [9, Eq. (2.1.6) and
(2.1.10)]). Additionally, Al-Salam and Carlitz provide formulas for the
moments of polynomials bearing their names, of which the discrete $q$%
-Hermite polynomials are special cases $\left[ 1\right] $. Therefore, it is
possible to give explicit representations of these polynomials as
determinants of moments. However, the fact that in [9, Eq. (2.1.10)] those
are Hankel determinants (or determinants close to Hankel determinants in [9,
(2.1.6]) clearly shows the difference from our formulas. For example, for $%
n=4$, we have$:$%
\begin{equation*}
H_{4}(x;q)=\left| 
\begin{tabular}{llll}
$x$ & $1-q$ & $0$ & $0$ \\ 
$1$ & $x$ & $1-q^{2}$ & $0$ \\ 
$x$ & $1$ & $x$ & $1-q^{3}$ \\ 
$1$ & $x$ & $1$ & $x$%
\end{tabular}
\right| \text{ \ \ \ by\ }\left( 9.25\right) \text{,}
\end{equation*}
while [9, Eq. (2.1.6)] gives after simplification: 
\begin{equation*}
H_{4}(x;q)=\dfrac{1}{q^{4}\left( 1-q^{2}\right) ^{2}}\left| 
\begin{tabular}{lllll}
$1$ & $0$ & $1-q$ & $0$ & $\left( 1-q\right) \left( 1-q^{3}\right) $ \\ 
$0$ & $1$ & $0$ & $\left( 1-q^{3}\right) $ & $0$ \\ 
$1$ & $0$ & $\left( 1-q^{3}\right) $ & $0$ & $\left( 1-q^{3}\right) \left(
1-q^{5}\right) $ \\ 
$0$ & $1$ & $0$ & $\left( 1-q^{5}\right) $ & $0$ \\ 
$1$ & $x$ & $x^{2}$ & $x^{3}$ & $x^{4}$%
\end{tabular}
\right|
\end{equation*}
\end{remark}

Consider now the generating function of the discrete q-Hermite II
polynomials denoted as $\widetilde{H_{n}}\left( x,q\right) $, [11, Eq.
(14.29.12)]: 
\begin{equation}
\widetilde{G}(t,x)=\sum\limits_{n\geq 0}\dfrac{\widetilde{H_{n}}\left(
x;q\right) }{\left( q;q\right) _{n}}q^{\binom{n}{2}}t^{n}=\dfrac{\left(
-xt;q\right) _{\infty }}{\left( -t^{2};q^{2}\right) _{\infty }}.  \tag{9.26}
\end{equation}
One can find determinant formulas for $\widetilde{H_{n}}$ using the formulas 
$\left( 9.24\right) $, $\left( 9.25\right) $ and $\left( 1.6\right) $, by
exploiting the fact that $H_{n}\left( x;q\right) =i^{-n}H_{n}\left(
ix;q^{-1}\right) $ (See $\left[ 11\text{, p. 552}\right] $). However, one
can also use the relations with complete symmetric functions. Indeed,
comparing $\left( 9.26\right) $ and $\left( 9.17\right) $ results in: 
\begin{equation}
\widetilde{G}(t,x)=\left( G\left( it,ix\right) \right) ^{-1}\text{.} 
\tag{9.27}
\end{equation}
Let $\Phi $ be the homomorphism corresponding to the specialization $G\left(
t,x\right) $ of $E(t)$. And let $\Phi \left( h_{n}\right) $ and $\Phi \left(
H\right) $ represent the images of the complete symmetric functions $h_{n}$
and the image of their generating function $H$, respectively. Using $\left(
2.4\right) $ and $\left( 9.27\right) $ we obtain: 
\begin{equation*}
\sum\limits_{n\geq 0}\dfrac{\widetilde{H_{n}}\left( x;q\right) }{\left(
q;q\right) _{n}}q^{\binom{n}{2}}t^{n}=\Phi \left( H\right) \left(
-it,ix\right) =\sum\limits_{n\geq 0}\left( \Phi (h_{n})\left( ix;q\right)
\right) \left( -it\right) ^{n}\text{,}
\end{equation*}
which implies: 
\begin{equation*}
\widetilde{H_{n}}\left( x;q\right) =\left( i(q-1)\right) ^{n}q^{-\binom{n}{2}%
}\left[ n\right] _{q}!\Phi (h_{n})\left( ix;q\right) \text{.}
\end{equation*}
Using Eq. $\left( 3.9\right) $, whith $\left[ p_{n}\right] $ given by $%
\left( 9.23\right) $ with $ix$ instead of $x$, one would obtain $\widetilde{%
H_{n}}(x;q)$ in determinant forms. Finally, by utilizing Eq. $\left(
6.11\right) $, we find the following $q$-exponential form of the generating
function:

\begin{equation*}
\sum\limits_{n\geq 0}\dfrac{\widetilde{H_{n}}\left( x;q\right) }{\left(
q;q\right) _{n}}q^{\binom{n}{2}}t^{n}=\mathbf{E}_{q}\left[ \left( 1-q\right)
^{-1}(\dfrac{xt}{\left[ 1\right] }-\dfrac{t^{2}}{\left[ 2\right] }-\dfrac{%
xt^{3}}{\left[ 3\right] }+\dfrac{t^{4}}{\left[ 4\right] }+\dfrac{xt^{5}}{%
\left[ 5\right] }-...)\right] _{q}^{\ast }\text{.}
\end{equation*}

\section{Additional developments}

In conclusion, here are other developments that have been realized or are
yet to be explored, in addition to those already mentionned:

*By applying the method from subsection 9.3 to Al-Salam-Carlitz polynomials,
determinant formulas and the $q$-exponential generating function have been
obtained for these polynomials. More generally, this method can be attempted
a priori for all $q$-othogonal polynomials with a generating function in the
form of an infinite product, as listed in [11].

* Similarly to what has been done in $\left[ 4\right] $, one can study the $%
p,q$-power symmetric functions $\left[ p_{n}\right] _{p,q}$ and $\left[
p_{\lambda }\right] _{p,q}$. It is easy to verify that many of the
properties established in the present article generalize to the $p,q$%
-analog. However, for certain properties such as factorizations into
infinite product, this does not seem possible.

* It should be investigated whether the relationships between $p_{\lambda }$ 
$\ $and the Schur function $s_{\lambda }$ have interesting $q$-analogs with $%
\left[ p_{\lambda }\right] $.

\bigskip

\textbf{References}

$\left[ 1\right] $ W. A. Al-Salam, L. Carlitz, Some Orthogonal $q$%
-Polynomials. \textit{Math. Nachr.}, \textbf{30}, (1965) 47-61.

$\left[ 2\right] $ F. Bergeron, \textit{Algebraic Combinatorics and
Coinvariant Spaces}, CRC Press, Ottawa, 2019.

$\left[ 3\right] $ Bourbaki, \textit{El\'{e}ments de math\'{e}matique,
Alg\`{e}bre}, Chap. 4, Springer, 1981.

$\left[ 4\right] $ V. Brugidou, A $q$-analog of certain symmetric functions
and one of its specializations, (2023); arXiv:2302.11221.

$\left[ 5\right] $ V. Brugidou, On a particular specialization of monomial
symmetric functions, (2023); arXiv:2306.15300.

$\left[ 6\right] $ L. Comtet, \textit{Advanced Combinatorics}. Springer
Netherlands, Dordrecht,1974.

$\left[ 7\right] $ G. Gasper, M. Rahman, \textit{Basic Hypergeometric series}%
, second ed., Cambridge University Press, 2004.

$\left[ 8\right] $ I. Gessel, A $q$-analog of the exponential formula. 
\textit{Discrete Mathematics}, \textbf{40} (1982) 69-79.

$\left[ 9\right] $ M. E. Ismail, \textit{Classical and Quantum Orthogonal
Polynomials in One Variable, }Cambridge University Press, 2005.

$\left[ 10\right] $ V. Kac, P. Cheung, \textit{Quantum Calculus. }Springer,
New York, 2002.

$\left[ 11\right] $ R.Koekoek, P.A Lesky, R.F. Swarttouw, \textit{%
Hypergeometric Orthogonal Polynomials and Their }$q$\textit{-Analogues},
Springer, 2010.

$\left[ 12\right] $ P. S. Kung, C. H. Yan, Goncarov polynomials and parking
functions. \textit{J. Comb. Theory Ser. A}, \textbf{102} (1), (2003), 16-37.

$\left[ 13\right] $ I. G. Macdonald, \textit{Symmetric functions and Hall
polynomials}, second ed., Oxford University Press, 1995.

$\left[ 14\right] $ C. L. Mallows and J. Riordan, The inversion enumerator
for labeled trees. \textit{Bull. Amer. Soc. }\textbf{74} (1968) 92-94.

$\left[ 15\right] $ R. P. Stanley, \textit{Enumerative Combinatorics.} Vol.
2, Cambridge University Press, 1999.

$\left[ 16\right] $ C.H. Yan, Parking functions, in M. Bona (ed.), \textit{%
Handbook of Enumerative Combinatorics}, CRC Press, Boca Raton, FL (2015) pp.
835-893.

\end{document}